\documentclass[runningheads]{llncs}
\usepackage{blindtext}
\usepackage{chemformula}
\usepackage{amssymb}
\usepackage{ifthen}
\usepackage{amsmath,times}
\usepackage{algorithm}
\usepackage{algorithmicx}
\usepackage{algpseudocode}
\usepackage{setspace}
\usepackage[percent]{overpic}
\usepackage{color}
\usepackage{graphicx,subfigure}
\usepackage{algpseudocode}
\usepackage{tikz}
\usepackage{hyperref}
 \usepackage{float} 
\begin{document}
\title{Conditioning by Projection for the Sampling from Prior Gaussian Distributions
}
%
%
\author{Alsadig Ali\inst{1}\and Abdullah Al-Mamun\inst{2}\and Felipe Pereira\inst{1} \and Arunasalam Rahunanthan\inst{3}}
\authorrunning{A. Ali et al.}
%
\institute{Department of Mathematical Sciences, The University of Texas at Dallas, Richardson, TX, 75080, USA \\ \and
Institute of Natural Sciences, United International University, Dhaka, Bangladesh\\ \and 
Department of Mathematics and Computer Science, Central State University,
Wilberforce, OH 45384, USA\\
\email{aRahunanthan@centralstate.edu}}
\maketitle              
%
\begin{abstract}

In this work we are interested in the (ill-posed) inverse problem for absolute permeability characterization 
that arises in predictive modeling of porous media flows. We consider a Bayesian statistical framework with
a preconditioned Markov Chain Monte Carlo (MCMC) algorithm for the solution of the inverse problem. Reduction of uncertainty
can be accomplished by incorporating measurements at sparse locations (static data) in the prior distribution. We present a new method to condition Gaussian fields (the log of permeability fields) to available sparse measurements. 
A truncated Karhunen-Lo\`eve expansion (KLE) is used for dimension reduction. In the proposed method the imposition of
static data is made through the projection of a sample (expressed as a vector of independent, identically distributed
normal random variables) onto the nullspace of a data matrix, that is defined in terms of the KLE. 
The numerical implementation of the proposed method is straightforward.
Through numerical experiments for a model of second-order elliptic equation, we show that the proposed method in multi-chain studies converges much faster than the MCMC method without conditioning. These studies indicate the importance of conditioning in accelerating the MCMC convergence.
\keywords{Conditioning  \and Predictive simulations \and Bayesian framework \and MCMC.}
\end{abstract}
%
%
\section{Introduction}

Subsurface formations have spatial variability in their hydraulic properties in multiple length scales. It 
has been established that such variability plays an important role in determining fluid flow patterns in heterogeneous 
formations \cite{glimm92,glimm93,furtado03,mrborges08}. The characterization of such formations is essential for 
predictive simulations of multiphase flows in problems such as oil recovery, \ch{CO2} sequestration, or contaminant transport. 
Procedures for characterization, such as a Bayesian framework using the Markov chain Monte Carlo (MCMC) algorithm, entails 
the solution of an ill-posed inverse problem. From static and dynamic data, the unknown coefficients of a governing system 
of partial differential equations that model fluid flow have to be recovered. 
In the inversion process, the available static data (for instance, measured absolute permeability values at sparse locations) 
have to be incorporated in the sampling of the (unknown) permeability field so that uncertainty is reduced as much as possible. 
This way,  a conditioning procedure has to be applied such that all samples considered by the inversion method honor 
available measurements.

MCMC methods have been applied in several areas of science and engineering, and have attracted the attention of
many research groups~\cite{Ginting2012707,abdullah2018}. There are many recent and important developments related to gradient-based MCMC procedures.
We mention, among many others, the Hamiltonian Monte Carlo \cite{HMC_2011}, and the active subspace methods \cite{Colorado_2016}. There are also 
Hessian-based procedures \cite{Omar_2012}. However, for subsurface flow problems gradient and Hessian calculations are computationally very 
expensive. Thus, for the problems we have in mind we consider the preconditioned 
(or two-stage) MCMC introduced 
in \cite{efendiev2006,christenfox05}. Some developments of this method include the multi-physics version \cite{Pereira_Rahu_2015}, the parallelization in multi-core devices \cite{Pre_fetching_Pereira_Rahu_2014} and its extension to deformable subsurface formations \cite{Pereira_Borges_2021}. This procedure has also been successfully applied in geophysics \cite{Georgia_2019}.
The Bayesian framework for inversion that we will use in this work has been carefully discussed in \cite{Abdullah_2020}. 

Our focus in this work is the conditioning of samples to sparse measured data. The importance of conditioning in reservoir 
simulation has been recognized for many years (see, for instance, \cite{Hewett_Behrens}).
The sequential Gaussian simulation (SGS) \cite{SGS_1992} has been widely used for this purpose. More recently a truncated
Karhunen-Lo\`eve expansion (KLE) \cite{loeve1977,efendiev2006,ginting2011_2} has been applied both for dimensional reduction of the 
(computationally infeasible) large stochastic dimensions of fine grid discretizations of subsurface flow problems as 
well as for conditioning. In \cite{Malgo} one can find a discussion of the advantages of the KLE over the SGS. 
There are two lines of work for the conditioning of samples using the KLE. In \cite{KLE_MOD_2004} the basic
KLE construction is altered to incorporate static data. In another line of work the vector of independent, identically distributed (i.i.d.) $\mathcal N(0,1)$ 
Gaussian variables that enter the KLE construction is modified to take data into account while keeping the KLE eigenvectors unaltered. This approach can be seen 
in \cite{mondal10} where a search through a set of
matrices of size determined by the number of measurements is performed, and they take the one with the best condition number.
In \cite{Malgo} the authors use the expressions for conditional means and covariances of Gaussian random variables \cite{ref-19-Malgo_2012}
that require the inversion of a data covariance matrix. 

The method that we introduce here does not alter the basic KLE construction and is computationally very inexpensive.
By using the measured data in the KLE we define a data matrix that models Gaussian perturbations on top of a surface
defined by kriging \cite{kriging_1978} of the available data. The proposed procedure consists of replacing the vectors of 
i.i.d., $\mathcal N(0,1)$ Gaussian variables by the closest vector to them (in the least-squares sense) in the nullspace of the data 
matrix. This way, the Gaussian fields produced by the KLE honor exactly all the data available. We remark that all one 
needs to implement the proposed method is a basis for the nullspace of the data matrix. Such a basis is used to construct 
the projection onto the nullspace of the data matrix. Besides introducing this new method, we also present multi-chain MCMC studies 
for a model elliptic equation. We show convergence of this method for large dimensional problems. Moreover, these studies
are used to illustrate the importance of conditioning for the convergence of MCMC methods. Although outside the scope of
this work we believe the method presented here can also be applied to other sampling strategies such as the hierarchical 
method discussed in \cite{Borges_Pereira_Hierarchical_2010}. This topic is currently being considered by the authors.

This work is organized as follows. In section~\ref{model} we discuss the model problem, kriging and dimensional reduction. Section~\ref{proj} introduces the new method for conditioning by projection. In section~\ref{Bayesian} we explain the Bayesian framework for inversion. We describe the preconditioned MCMC method with and without conditioning for the construction of permeability fields. Also, the method we use for convergence assessment of MCMC methods is presented. Section~\ref{results} is dedicated to a discussion of numerical results from our experiments. Conclusions and some remarks appear in section~\ref{conclusions}.

\section{Model Problem and Dimensional Reduction}
\label{model}
\subsection{The Model Problem}

In the modeling of incompressible subsurface flow problems, such as in two-phase flows in oil reservoirs \cite{Douglas_Pereira_MCMCAA_1997,Rahu_Pereira_MATCOM_2011}, the system of governing partial differential equations consists of a second order elliptic equation coupled to a hyperbolic dominated
equation that describes fluid transport. In this work, we  illustrate the proposed method for conditioning Gaussian fields in terms of the elliptic equation that enters in this system.

Let us consider $\Omega \subset \mathbb{R}^2$, a bounded domain with Lipschitz boundary, let $\boldsymbol{v}(\boldsymbol{x})$ and $p(\boldsymbol{x})$ denote the Darcy velocity and the  pressure of the fluid, respectively, where $\boldsymbol{x} \in \Omega$, and $k(\boldsymbol{x})$ is the absolute permeability. Using Darcy's law, we have the following model  elliptic problem:
\begin{equation}
\begin{aligned}
\label{p_eqn}
\boldsymbol{v}(\boldsymbol{x})~~~~~ &= -k(\boldsymbol{x}) \nabla p(\boldsymbol{x}),\\
\nabla \cdot \boldsymbol{v}(\boldsymbol{x}) &= f,~~~ \boldsymbol{x} \in \Omega,
\end{aligned}
\end{equation}
where $f\in L^2(\Omega)$ is a given source term, together with appropriate Dirichlet and Neumann type boundary conditions $p= g_p\in H^{\frac{1}{2}}(\partial \Omega_p)$ and $\boldsymbol{v}\cdot\boldsymbol{\hat{n}}=g _{\boldsymbol{v}}\in H^{-\frac{1}{2}}(\partial \Omega_{\boldsymbol{v}})$, the pressure boundary data and normal velocity data, respectively, with $\partial\Omega = \overline{\partial\Omega_p} \cup \overline{\partial\Omega_{\boldsymbol{v}}}$ and $\overline{\partial\Omega_p} \cap  \overline{\partial\Omega_{\boldsymbol{v}}} = \varnothing$, $\boldsymbol{\hat{n}}$ is the exterior unit normal vector. The numerical approximation of the above elliptic problem is derived from the weak formulation of the velocity-pressure system. The weak formulation can be found in \cite{mumm_2014}.

\subsection{Kriging}
\label{kriging}
Kriging is a type of interpolation method, which uses a local estimation technique that gives the best linear unbiased estimator of unknown values~\cite{kriging_1978,mining_1976}. The main idea of kriging is to estimate the values of a function (or random field) in a domain by calculating a weighted average of a given point in some neighborhood. The kriging takes into account (i) the distances between the estimated points in the domain and the given data points, and (ii) the distances between the given data points themselves. Thus, the kriging gives more weights to the nearest data points. When a point $\boldsymbol{x}_0$ coincides with a given data point $\boldsymbol{x}_i$, the kriging gives the exact value and the variance is zero. In the following subsection, we discuss the KLE, and in Section \ref{proj}, we combine both kriging and KLE.

\subsection{Dimensional Reduction}
We use the KLE to reduce the large dimensional uncertainty space describing the permeability field $k(\boldsymbol{x})$. In this subsection, we present a brief description of the KLE \cite{loeve1977,ginting2013}.
Let $\boldsymbol{x}\in \Omega$, and suppose $\text{log}\left[ k(\boldsymbol{x})\right]= Y^k(\boldsymbol{x})$ is a Gaussian field. Also, let us assume $Y^k(\boldsymbol{x})$ is a second-order stochastic process, that is $Y^k(\boldsymbol{x})\in L^2(\Omega)$ with a probability one. Set $E[(Y^k)^2] = 0$. Thus, for a given orthonormal basis $\left \{\varphi_i\right\}$ in $L^2(\Omega)$, $Y^k(\boldsymbol{x})$ can be written as the following:
\begin{equation}
\label{kle1}
Y^k(\boldsymbol{x}) = \sum _{i=1}^{\infty} Y_i^k \varphi_i(\boldsymbol{x}), \hspace*{0.3cm}
\text{with random coefficients }\hspace*{0.3cm}
Y_i^k = \int_{\Omega}Y^k(\boldsymbol{x}) \varphi_i(\boldsymbol{x})d\boldsymbol{x},
\end{equation}
where $\varphi_i(\boldsymbol{x})$ is an eigenfunction satisfying 
\begin{equation}
\label{kle_efun}
\int_{\Omega}R(\boldsymbol{x}_1, \boldsymbol{x}_2)\varphi_i(\boldsymbol{x}_2)d\boldsymbol{x}_2 = \lambda_i \varphi_i(\boldsymbol{x}_1),~~ i = 1, 2, . . .,
\end{equation}
with $\lambda_i = E[(Y_i^k)^2] > 0$ and $R(\boldsymbol{x}_1,\boldsymbol{x}_2)$ is the covariance function, Cov$(Y^k(\boldsymbol{x}_1),Y^k(\boldsymbol{x}_2))$. Letting $\theta_i^k = Y_i^k/\sqrt{\lambda_i}$, the KLE in equation \eqref{kle1} can be written as  
\begin{equation}
\label{kle_exp}
Y^k(\boldsymbol{x}) = \sum_{i=1}^{\infty} \sqrt{\lambda_i}{\theta}_i^k\varphi_i (\boldsymbol{x}),
\end{equation}
where $\lambda_i$ and $\varphi_i$ satisfy equation \eqref{kle_efun}.  
We arrange the eigenvalues in the descending order. The infinite series in equation \eqref{kle_exp} is called the Karhunen-Lo\`eve expansion. In general, we only need first $n$ dominating terms such that the energy $E$ is more than $95\%$~\cite{Laloy2014}, i.e. 
\begin{equation}
\label{energy}
E=\dfrac{\sum_{i=1}^n \lambda_i}{\sum_{i=1}^{\infty} \lambda_i}\geq 95\% .
\end{equation}
Then, a truncated KLE is defined by
\begin{equation}
\label{kle_truncated}
Y^k_{n}(\boldsymbol{x})=\sum_{i=1}^{n}\sqrt{\lambda_i}\theta_i^k
\varphi_i(\boldsymbol{x}).
\end{equation}  
From now on we will use $Y(\boldsymbol{x})$ instead of $Y_n^k(\boldsymbol{x})$ for samples generated by the truncated KLE.

\section{Conditioning by Projection}
\label{proj}
Suppose the Gaussian field $Y(\boldsymbol{x})$ defined in equation~\eqref{kle_truncated} is a Gaussian perturbation on top of a kriged field $\hat{Y}(\boldsymbol{x})$ discussed in section~\ref{kriging}. That is,
\begin{equation}
\begin{aligned}
\label{condi_perm}
Y(\boldsymbol{x}) = \hat{Y}(\boldsymbol{x})+ \sum_{i=1}^{n}\sqrt{\lambda_i}
\varphi_i(\boldsymbol{x})\theta_i ~.
\end{aligned}  
\end{equation}
This can also be written as
\begin{equation}
\begin{aligned}
\label{condi_perm1}
Y(\boldsymbol{x}) - \hat{Y}(\boldsymbol{x}) &= \sum_{i=1}^{n}\sqrt{\lambda_i}\varphi_i(\boldsymbol{x})\theta_i\\
  &=\pmb{\phi}^T(\boldsymbol{x})\sqrt{D}\pmb{\theta},
\end{aligned}  
\end{equation}
where $D$ is a diagonal matrix of the dominating $n$ eigenvalues and $\pmb{\theta}=(\theta_1,\theta_2,\dots, \theta_n)$. For each $\boldsymbol{x}$, $\pmb{\phi}(\boldsymbol{x})\in \mathbb{R}^n$.
Furthermore, assume that $m$ measured data values $Y^1,Y^2,\dots ,Y^m$ of $Y(\boldsymbol{x})$ are given at sparse locations $\boldsymbol{\hat{x}}=(\boldsymbol{\hat{x}}_1,\boldsymbol{\hat{x}}_2,\dots ,\boldsymbol{\hat{x}}_m)$ in the domain $\Omega$ such that 
\begin{equation}
    \label{measured_data}
    Y(\boldsymbol{\hat{x}}_i)=Y^i,\quad  i=1,2,\dots,m
\end{equation}
are known. Thus, by substituting~\eqref{measured_data} in equation~\eqref{condi_perm1} we have 
\begin{equation}
\begin{aligned}
\label{condi_perm2}
\pmb{\phi}^T(\boldsymbol{\hat{x}})\sqrt{D}\pmb{\theta} &= Y(\boldsymbol{\hat{x}}) - \hat{Y}(\boldsymbol{\hat{x}}) \\
  &= \boldsymbol{0}.
\end{aligned}  
\end{equation}
This honors $m$ data values at the corresponding sparse locations $\boldsymbol{\hat{x}}$. We can rewrite~\eqref{condi_perm2} as follows
\begin{equation}
\begin{aligned}
\label{system}
A\pmb{\theta}= \boldsymbol{0}
\end{aligned}  
\end{equation}
where $A=\pmb{\phi}^T(\boldsymbol{\hat{x}})\sqrt{D}\in \mathbb{R}^{m\times n}$ is our data matrix. The solution of this homogeneous system of equations~\eqref{system} is the nullspace $N(A)$, a subspace of $\mathbb{R}^n$. Given an arbitrary vector $\pmb{\theta}\in \mathbb{R}^n$ of i.i.d., $\mathcal N(0,1)$ Gaussian variables that would be used in the KLE~\eqref{condi_perm}, in order to honor the given $m$ data we select the closest vector to $\pmb{\theta}$ in the subspace $N(A)$. This vector is obtained by projecting $\pmb{\theta}$ onto the nullspace of the data matrix $A$. 

Assume that $A$ has rank $r$, we can easily find an orthonormal basis $\beta = \{\boldsymbol{q}_1,\boldsymbol{q}_2,\dots,\\ \boldsymbol{q}_{n-r}\}$ for the nullspace of $A$. We assume that $r<n$, so that $A$ has a non-trivial nullspace. The projection matrix $P$ onto the nullspace of $A$ is given by~\cite{LA_strang_2019}
\begin{equation}
P = QQ^T,
\end{equation}
where $Q$ is a matrix that has $\boldsymbol{q}_1,\boldsymbol{q}_2,\dots,\boldsymbol{q}_{n-r}$ as columns. The projection vector of $\pmb{\theta}$ onto the nullspace of the data matrix $A$ is given by 
\begin{equation}
\label{projection}
\pmb{\hat{\theta}} = P\pmb{\theta}.
\end{equation}
\vspace{-1.4cm}
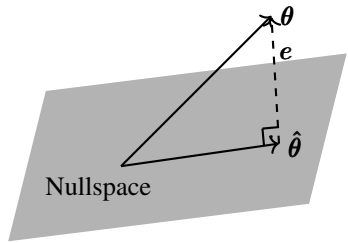
\begin{figure}[H]
	\centering
	\begin{tikzpicture}[scale =1]
	\fill[black!30!white] (-1,1)--(3,1.5)--(2.5,-0.5)--(-1.5,-1)--(-1,1);
	\draw[thick][->](0,0)--(2,2) node[right]{$\pmb{\theta}$};
	\draw[thick][->](0,0)--(2.1,0.3) node[right]{$\pmb{\hat{\theta}}$};
	\draw[thick][dashed][->] (2.1,0.3)--(2,1.98);
	\node at (2.2,1.5) {$\pmb{e}$};
	\node at (-0.3,-0.3) {Nullspace};
	\draw[thick][-](1.9,0.28)--(1.88,0.5);
	\draw[thick][-] (1.88,0.5)--(2.1,0.52);
	\node at (2,0.4) {$.$};
	\end{tikzpicture}
	\caption{Projection of the i.i.d., $\mathcal N(0,1)\,\pmb{\theta} $ showing the  closet vector $\pmb{\hat{\theta}}$ in the nullspace.}
	\label{null_space}
\end{figure}
Fig.~\ref{null_space} illustrates the projection of the i.i.d., $\mathcal N(0,1)$ Gaussian variables $\pmb{\theta}$ vector onto the nullspace of the data matrix $A$. The error vector defined as $\pmb{e}=\pmb{\theta}-\pmb{\hat{\theta}}$ is minimized in the least squares sense. It is perpendicular to the projection vector $\pmb{\hat{\theta}}$. Finally, we can rewrite equation~\eqref{condi_perm} as  
\begin{equation}
\begin{aligned}
\label{condi_perm3}
Y(\boldsymbol{x}) = \hat{Y}(\boldsymbol{x})+ \sum_{i=1}^{n}\sqrt{\lambda_i}
\varphi_i(\boldsymbol{x})\hat{\theta_i},
\end{aligned}  
\end{equation}
where $\hat{\theta_i},i=1, \dots ,n$ are the components of $\pmb{\hat{\theta}}$ in equation~\eqref{projection}. If $\boldsymbol{x}=\boldsymbol{\hat{x}}_i, i=1,2,\dots, m$, in equation~\eqref{condi_perm3}, the second term on the right-hand side is zero. This leads to $Y(\boldsymbol{x}_i)=\hat{Y}(\boldsymbol{\hat{x}}_i)$, and thus, we honor the exact values of the permeability field at the known locations $\boldsymbol{\hat{x}}_1,\dots, \boldsymbol{\hat{x}}_m$.

\section{The Bayesian Framework}
\label{Bayesian}

\subsection{Posterior Exploration}
\label{post}
Our primary interest in this article is to characterize the permeability field of our domain of study. Here, we discuss the characterization of the permeability field conditioned on available pressure data. We use half (in a chessboard pattern) of the total pressure data in the discretized domain~\cite{mala_2020}. We denote the (log of the) permeability field and the reference pressure data by $\pmb {\eta}$ and $R_p$, respectively. A Bayesian statistical approach along with a preconditioned MCMC method with and without conditioning is used to solve the inverse problem. The posterior probability is calculated using the Bayes' theorem with respect to $R_p$ as follows:  
\begin{equation}
  \label{bayes_eqn}
 P(\pmb{\eta}|R_p) \propto P(R_p|\pmb{\eta})P(\pmb{\eta}),
\end{equation} 
where $P(\pmb{\eta})$ represents the Gaussian prior distribution. The normalizing constant is disregarded in this investigation as we do an iterative search in our MCMC method. We construct the permeability field $\pmb{\eta}(\pmb{\theta})$ using the KLE, where $\pmb{\theta}$ is generated by the MCMC method. The likelihood is considered to be a Gaussian function as in~\cite{efendiev2006}, and as a result of that we have our likelihood function as
\begin{equation}
  \label{likelihood_fun}
  P(R_p|\pmb{\eta}) \propto \exp\Big(-(R_p - R_{\pmb\eta})^\top\Sigma (R_p - R_{\pmb\eta})\Big),
\end{equation}
where $R_{\pmb\eta}$ denotes the simulated pressure data. We set the covariance matrix  $\Sigma$  as $\Sigma =  \pmb{I}/2\sigma_R^2$, where $\pmb{I}$ and $\sigma_R^2$ represents the identity matrix and the precision parameter, respectively.

We use a MCMC algorithm to sample from the posterior distribution~\eqref{bayes_eqn}.  In the MCMC algorithm, we have an instrumental distribution $I(\pmb{\eta}_p|\pmb{\eta})$ by which we propose $\pmb{\eta}_p =\pmb{\eta}(\pmb{\theta}_p)$ at every iteration, where $\pmb{\eta}$ represents the previously accepted state. We solve the forward problem for a given permeability field on the numerical simulator and compute the acceptance probability of a proposed proposal using   
\begin{equation}
	\label{single_stage_prob}
  {\alpha}(\pmb{\eta}, \pmb{\eta}_p) = \text{min}
  \left(1,\frac{I(\pmb{\eta}|\pmb{\eta}_p)P(\pmb{\eta}_p|R_p)}{I(\pmb{\eta}_p|\pmb{\eta})P(\pmb{\eta}|R_p)}\right).
\end{equation} 
Below we describe the preconditioned MCMC method with and without conditioning. 

\subsection{MCMC algorithms with and without conditioning}
\label{two-stage}
In our study, we use the MCMC method with and without conditioning. The following random walk sampler (RWS)~\cite{Cotter_2013} is used for both MCMC cases: 
\begin{equation}
\label{two_sampler}
\boldsymbol\theta_p = \sqrt{1-\beta^2}\, \boldsymbol\theta + \beta\,\boldsymbol\epsilon, 
\end{equation} 
where $\boldsymbol\theta_p$ denotes the current proposal and $\boldsymbol\theta$ represents the previously accepted proposal. The symbol $\beta$ denotes the tuning parameter and $\boldsymbol\epsilon $ represents $\mathcal N(0,1)$-random variable.
\begin{algorithm}
  \caption{MCMC with conditioning}
  \label{alg_two_stage}
   \begin{algorithmic}[1]
			\State For a given covariance function $R$ generate KLE.
				\For{$i=1$ to $M_{\text{mcmc}}$}
						
							\State At $\pmb{\eta}(\pmb{\theta})$generate $\pmb{\theta}_p$ using equation~\eqref{two_sampler}.
							\State Project $\pmb{\theta}_p$ to the nullspace using~\eqref{projection} to get $\pmb{\hat{\theta}}_p$.
							\State At $\pmb{\hat{\theta}}_p$ construct $\pmb{\eta}_p$ using equation~\eqref{condi_perm3}.
							\State Compute the upscaled permeability on the coarse-scale using $\pmb{\eta}_p$.
							\State Solve the forward problem on the coarse-scale to get $R_c$.
							\State Compute the coarse-scale acceptance probability ${\alpha}_c(\pmb{\eta}, \pmb{\eta}_p)$.
				
        \If{ $\pmb{\eta}_p$ is accepted}
						\State Use $\pmb{\eta}_p$ in the fine-scale simulation to get $R_f$.
            \State Compute the fine-scale acceptance probability ${\alpha}_f(\pmb{\eta},\pmb{\eta}_p)$.
						\If{ $\pmb{\eta}_p$ is accepted}
							  $\pmb{\eta} = \pmb{\eta}_p$. 
						\EndIf		
				\EndIf		
				
     \EndFor
 \end{algorithmic}
\end{algorithm}

We first discuss the algorithm of the MCMC method without conditioning~\cite{efendiev2006,christenfox05}. The filtering step of this method is based on a coarse-scale model approximation of the governing equation~\eqref{p_eqn}. The coarse-scale discretization is similar to the fine-scale discretization and the permeability field $\pmb{\eta}(\pmb{\theta})$ is projected on the coarse-scale. An upscaling procedure~\cite{durlofsky1991} is used to set an effective permeability field on a coarse-grid that provides a similar average response as that of the underlying fine-scale problem. The numerical simulator is run on the coarse-scale model and produces the coarse-grid pressure field $R_c$. The coarse-scale and fine-scale acceptance probabilities are estimated as \begin{equation}
\label{two_prob}
	\begin{aligned}
			{\alpha}_c(\pmb{\eta}, \pmb{\eta}_p) &= \text{min}\left(1,\dfrac{I(\pmb{\eta}|\pmb{\eta}_p)P_c(\pmb{\eta}_p|R_p)}{I(\pmb{\eta}_p|\pmb{\eta})P_c(\pmb{\eta}|R_p)}\right)\text{, and}\\
			{\alpha}_f(\pmb{\eta},\pmb{\eta}_p) &= \text{min}\left(1,\dfrac{P_f(\pmb{\eta}_p|R_p)P_c(\pmb{\eta}|R_p)}{P_f(\pmb{\eta}|R_p)P_c(\pmb{\eta}_p|R_p)}\right),
	\end{aligned}		
\end{equation}
where $P_c$ and $P_f$ represent the posterior probabilities computed at coarse- and fine-scale, respectively. In the MCMC method with conditioning, we construct the permeability field  $\pmb\eta(\pmb\theta)$ using equation \eqref{condi_perm3}. The algorithm of the MCMC method without conditioning is described in~\cite{ICCS_2020} and the algorithm of the MCMC method with conditioning is presented in Algorithm~\ref{alg_two_stage}.

\subsection{Convergence Assessment}
\label{convergence_ass}
In this subsection, we discuss the convergence of the MCMC methods for a reliable characterization of the permeability field. Determining the stopping criterion in the MCMC simulations is a challenging task. For this reason, researchers use several MCMC convergence diagnostics to decide when it is safe to stop the simulations. Different diagnostic measures are discussed in \cite{Roy2020,cowles1996,mengersen1999,brooks1998}. A common approach is to start more than one MCMC chain using distinct initial conditions. In this paper, we study the convergence of the proposed MCMC methods using the Potential Scale Reduction Factor (PSRF) and its multivariate extension MPSRF ~\cite{Brooksgelman1998}. Note that the PSRF does not take all the parameters into account and thus may not provide a very reliable indicator of the convergence of the parallel chains. However, the MPSRF takes all the parameters and their interactions into consideration. Therefore, the MPSRF is a better measure to determine the convergence of several chains that explore a high dimensional parameter space.

Let us consider $l$ posterior draws of $\boldsymbol\theta$ (dimension $n$) for each of the $k>1$ parallel chains. Let  $\boldsymbol {\theta}_j^{c}$ be the generated samples at iteration $c$ in $j$th chain of $k$ parallel MCMCs. Then, the posterior variance-covariance matrix is computed by 
\begin{equation}
\label{variance_covarinance}
   \mathbf{\widehat{V}} = \frac{l-1}{l}\mathbf{W} + \left( 1+ \frac{1}{k}\right)\frac{\mathbf{B}}{l}.
\end{equation}
The within-chain covariance matrix $\mathbf{W}$ and between-chain covariance matrix $\mathbf{B}$ are defined by 
\begin{equation}
  \mathbf{W} = \frac{1}{k(l-1)} \sum\limits_{j=1}^k \sum\limits_{c=1}^l \left(\boldsymbol {\theta}_j^{c} - \boldsymbol {\bar \theta}_{j.}\right) \left(\boldsymbol {\theta}_j^{c} - \boldsymbol {\bar \theta}_{j.}\right) ^{T},
\end{equation}
and
\begin{equation}
    \mathbf{B} = \frac{l}{k-1} \sum\limits_{j=1}^k \left(\boldsymbol {\bar \theta}_{j.}-\boldsymbol {\bar \theta}_{..}\right) \left(\boldsymbol {\bar \theta}_{j.}-\boldsymbol {\bar \theta}_{..}\right)^{T},
\end{equation}
respectively. $\boldsymbol {\bar \theta}_{j.}$ and  $\boldsymbol {\bar \theta}_{..}$ represent the mean within the chain and the mean of the $k$ combined chains, respectively, and $T$ stands for transpose. The PSRFs are computed by
\begin{equation}
  \begin{aligned}
\text{PSRF}_\text{i} =\sqrt{\frac {\text{diag}(\mathbf{\widehat{V}})_i}{\text{diag}(\mathbf{W})_i}}, ~~~~\text{where} \, \,\, i=1,2,...,n.
 \end{aligned} 
 \end{equation}
 The maximum of PSRF values, which is closer to 1, indicates that MCMC proposals sample from approximately the same posterior distribution. The MPSRF is defined as~\cite{Malyshkina2008} 
\begin{equation}
\begin{aligned}
\text{MPSRF}&=\sqrt{ \left(\frac{l-1}{l} + \left(\frac{k+1}{k}\right) \lambda\right)},
 \end{aligned} 
 \end{equation}
where $\lambda$ is the largest eigenvalue of the positive definite matrix  $\mathbf{W}^{-1} \mathbf{B}/l$. If the proposals in $k$ chains are sampled from almost the same distribution, the between chain covariance matrix 
becomes zero. Consequently, $\lambda$ and the MPSRF approach to $0$, and $1$, respectively, for sufficiently large $l$. Thus, we achieve the required convergence of the chains.
\section{Numerical Results}
\label{results}
\subsection{Setup of the problem}
\begin{figure}[H]
	\centering
	\includegraphics[width= 2.3 in]{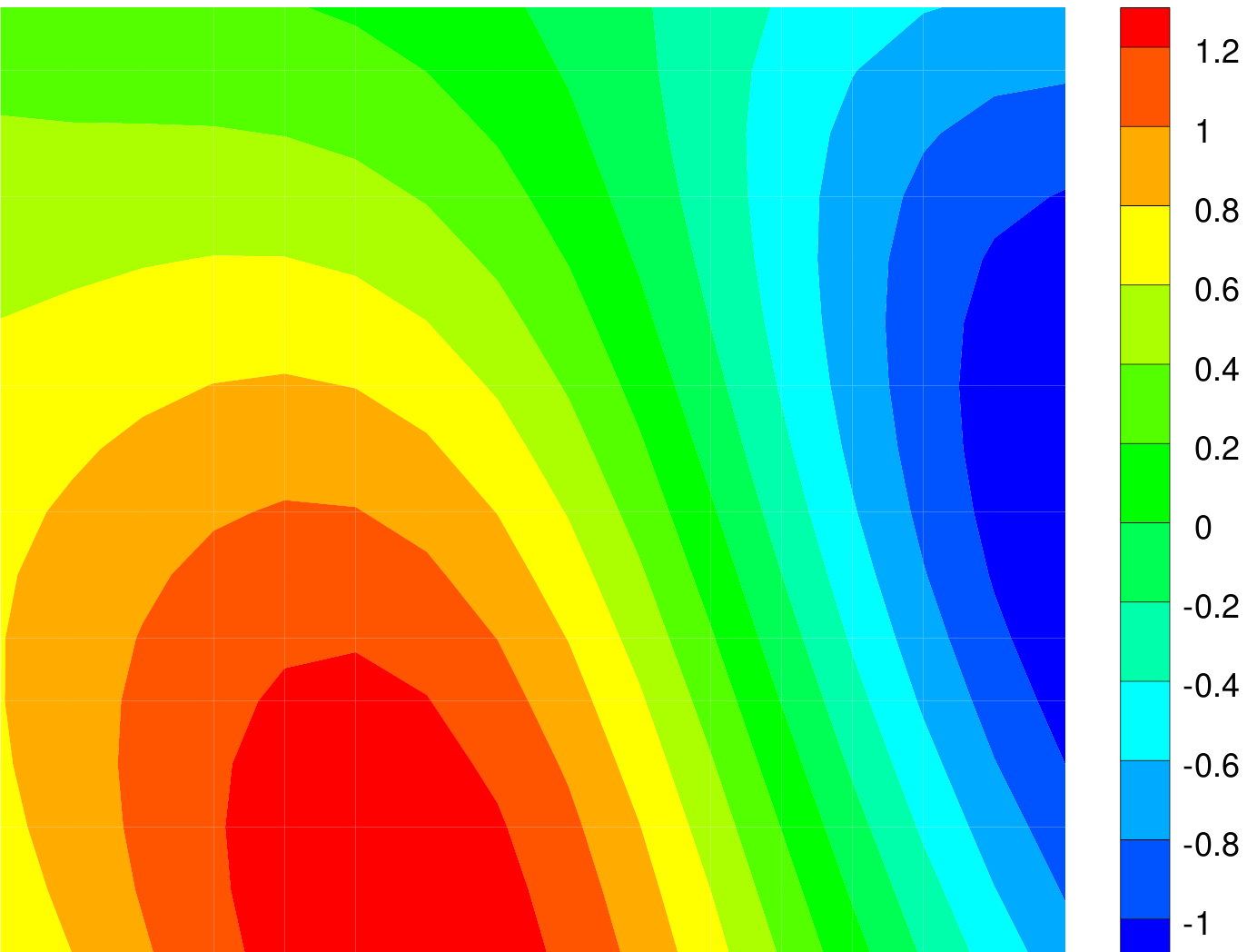}
	\hspace{2mm}
	\includegraphics[width= 2.3 in]{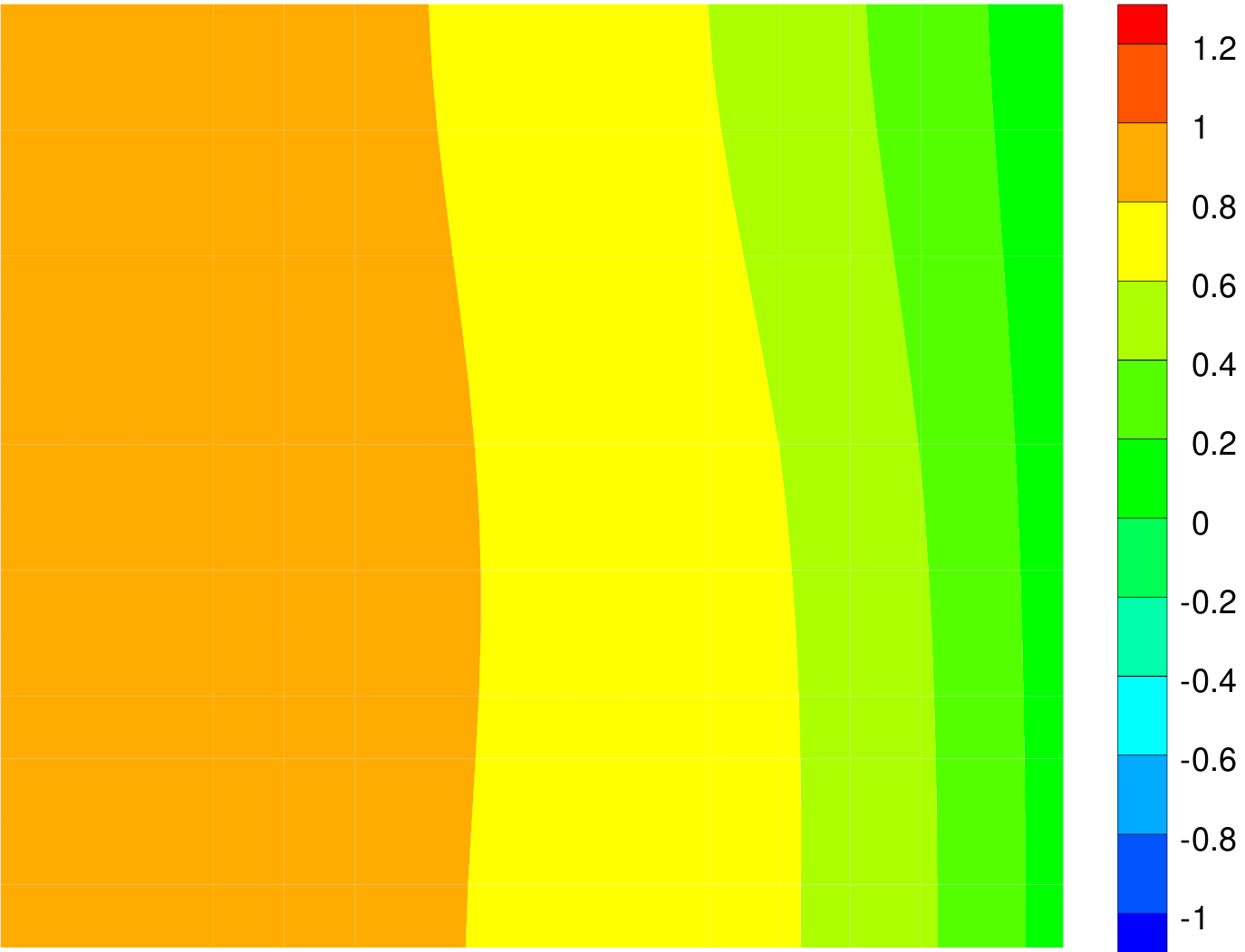}
	\caption{Reference log permeability field (left) and the corresponding reference pressure field (right).}
	\label{ref_perm}
\end{figure}

\begin{figure}[H] 
	\centering 
	\includegraphics[scale=0.8]{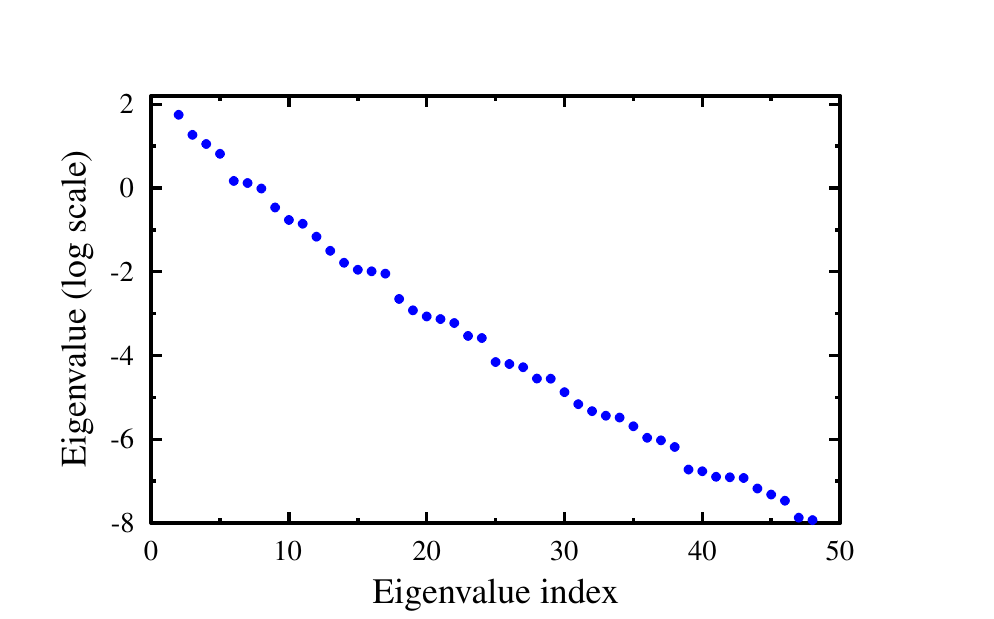}
	\caption{Decay of eigenvalues in KLE.}
	\label{eig_plot}
\end{figure}

In this subsection, we present the details of the simulation of our problem of interest. We consider a unit square-shaped physical domain and solve the elliptic equation~\eqref{p_eqn} in the domain. The simulated elliptic solutions are used in the MCMC algorithms. We then present a comparative study between the MCMC methods with and without conditioning. In the MCMC method with conditioning, we include the prior measurement of the permeability field at nine specific sparse locations. This permeability field is constructed using the KLE in which we consider the covariance function as
\begin{equation}
\label{kle}
\begin{aligned}
R(\boldsymbol{x}_1,\boldsymbol{x}_2) = \sigma_Y^2\, \text{exp}\left(-\frac{|x_1 - x_2|^2}{2L_x^2} - \frac{|y_1 - y_2|^2}{2 L_y^2}\right),
\end{aligned}
\end{equation}
where $L_x = 0.4$ and $L_y = 0.8$ are the correlation lengths and $\sigma_Y^2= \text{Var}[(Y^k)^2] = 1$. Fig.~\ref{eig_plot} illustrates the decay of the eigenvalues for these values. Considering the fast decay of the eigenvalues, in this paper, we take the first twenty dominating eigenvalues in the KLE, which preserves almost 100\% (using the equation~\eqref{energy}) of the total variance or energy of $Y$ in the equation~\eqref{kle_truncated}.
We set the source term $f=0$ and impose Dirichlet boundary conditions, $p=1$ and $p=0$, on the left and right boundaries, respectively. We set a Neumann-type boundary condition to zero (i.e., no-flow) everywhere else on the domain boundaries. In our study, we use a synthetic reference permeability field that is generated on a fine-grid of size $16\times16$. We then run the forward problem using the numerical simulator and generate the corresponding reference pressure field. See Fig.~\ref{ref_perm} for both fields.
We take the tuning parameter $\beta = 0.85$ in~\eqref{two_sampler} for both MCMC methods. As a screening step in the MCMC algorithm, we introduce a coarse-scale filter of size $8\times8$. The coarse-scale simulation runs almost two times faster than the corresponding fine-scale simulation while maintaining the general trend of the fine-scale simulation. We change one stochastic parameter at every iteration to construct a new permeability field for both MCMC studies.
\subsection{ Convergence Analysis of the MCMC Method}
\label{mpsrf_psrf}
We analyze here the convergence of the MCMC algorithm for both studies in terms of the maximum of PSRFs and MPSRF. We run four MCMCs starting from distinct initial seeds for the MCMC method without conditioning. We also run another four chains, but using the same corresponding initial seeds in the former runs, for the MCMC method with conditioning. In Table~\ref{rate} we show the acceptance rates for these two MCMC studies. We take $\sigma_F^2 = 10^{-4}$ and $\sigma_C^2 = 5\times 10^{-3}$ for coarse- and fine-scale precisions, respectively, in equation~\eqref{likelihood_fun}. 

\begin{figure}
  \centering
   \includegraphics[width = 2.36 in]{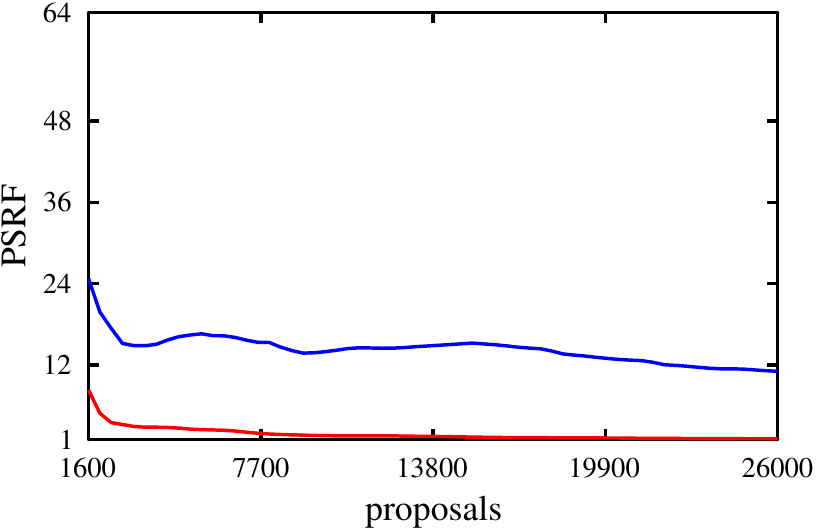}
 \hspace{0.1 mm}
   \includegraphics[width =2.36 in]{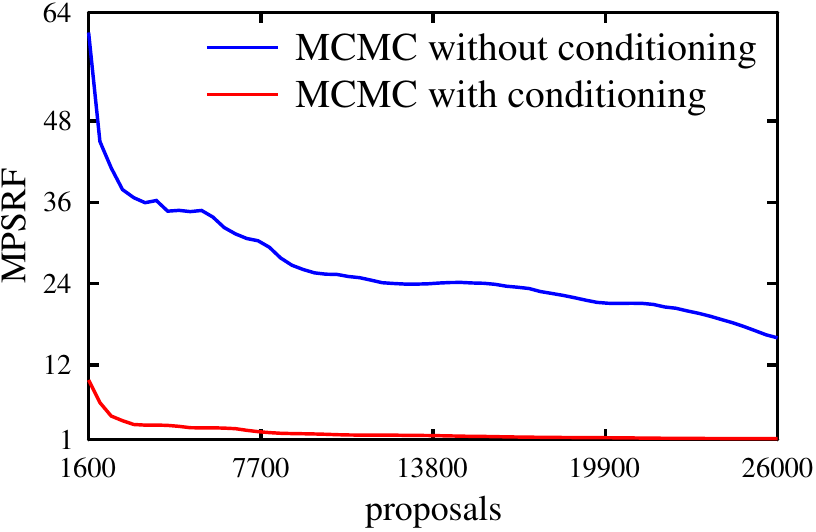}
   \caption{ The maximum of PSRFs and MPSRF for the MCMC method with and without  conditioning on measured data. }
   \label{MPSRF_PSRF_LS_2x2_20}
\end{figure}

\begin{table}[H] 
    \center
     \begin{tabular}{|cccc|}
	\hline
         &  \quad  MCMC without conditioning & \quad MCMC with conditioning&  \\
	\hline
         $\sigma_F^2$ &   $10^{-4}$  & $10^{-4}$  & \\               
         $\sigma_C^2$    & $5\times10^{-3}$ & $5\times10^{-3}$&   \\  
         acceptance rate                     & $53\%$   &$60\%$&  \\
        \hline
   \end{tabular}
\vspace{3 mm}
   \caption{A comparison of acceptance rates for the MCMC with and without conditioning.}
   \label{rate}     
\end{table}
We observe from the Table~\ref{rate} that the MCMC method with conditioning has a better acceptance rate than the MCMC method without conditioning.  In the calculation of the PSRFs and MPSRF, we take into account the $\pmb\theta$ values that are accepted on the fine-scale. If there is no change in the theta values from one state to another in the Markov chain, the repetition of the accepted theta values on the fine-scale is considered in the calculation. Note that a MCMC method converges to the target distribution if both PSRFs and MPSRF values are below $1.2$~\cite{Brian2007}. Fig.~\ref{MPSRF_PSRF_LS_2x2_20} depicts the maximum of the PSRFs and MPSRF for both cases with the same number of proposals, 26000, after burn-in. In each chain, 6500 proposals are considered. We observe that the MCMC with conditioning achieves a full convergence. Towards the end of the curves, the maximum PSRF and MPSRF values are 1.14 and  1.16, respectively. Both PSRF and MPSRF curves for the MCMC method without conditioning are considerably away from the numerical value of $1.2$. Thus, we conclude that the MCMC method with conditioning explores the posterior distribution much faster than the MCMC method without conditioning. Moreover, Fig.~\ref{MPSRF_PSRF_LS_2x2_20} demonstrates that the maximum of PSRFs is bounded above by the MPSRF as shown in~\cite{Brooksgelman1998}. In the next paragraph, we discuss the importance of this convergence in order to recover the underlying permeability field.

\begin{figure}
	\centering	
	\includegraphics[width=1.575in]{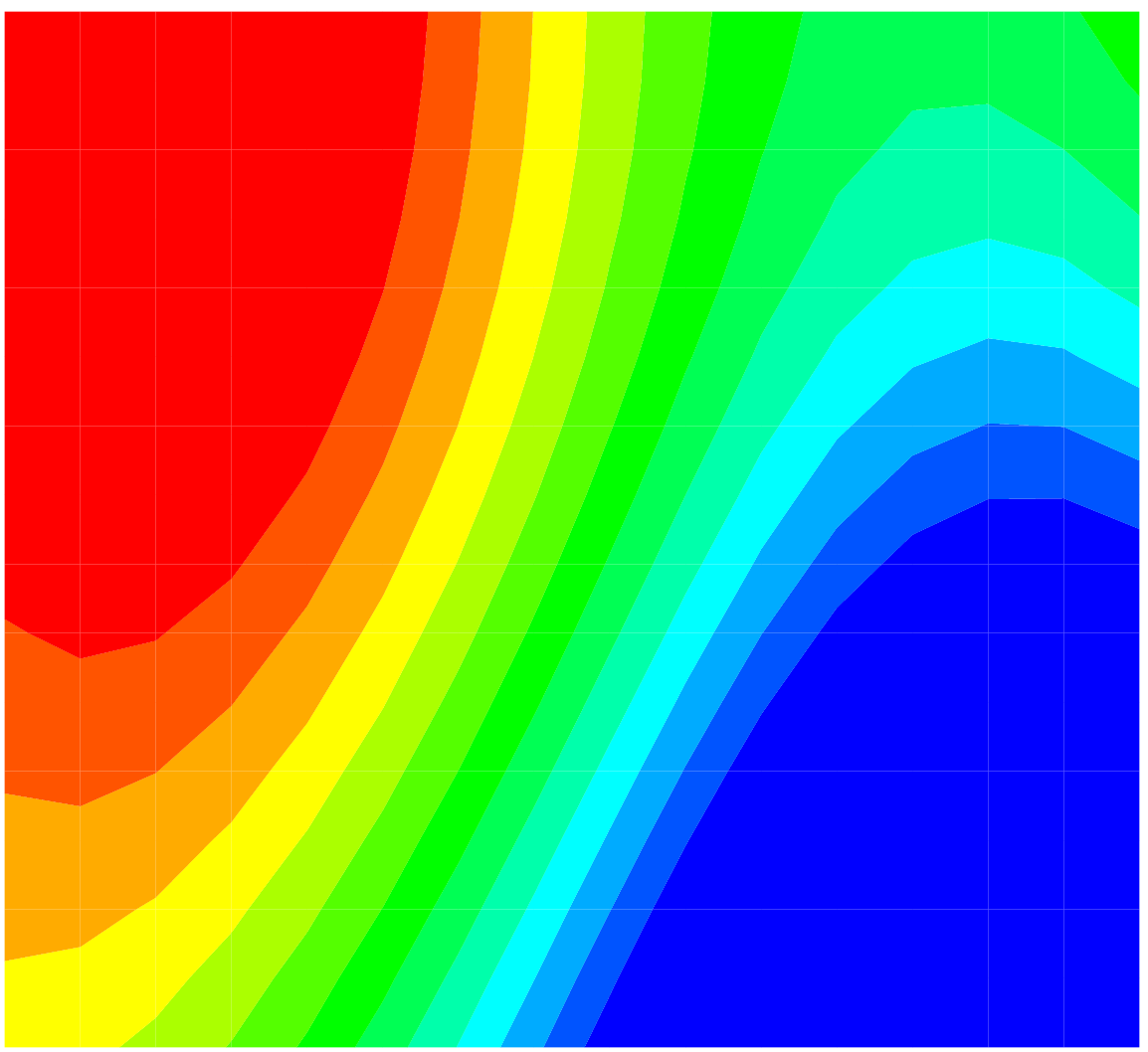}
	\includegraphics[width=1.575in]{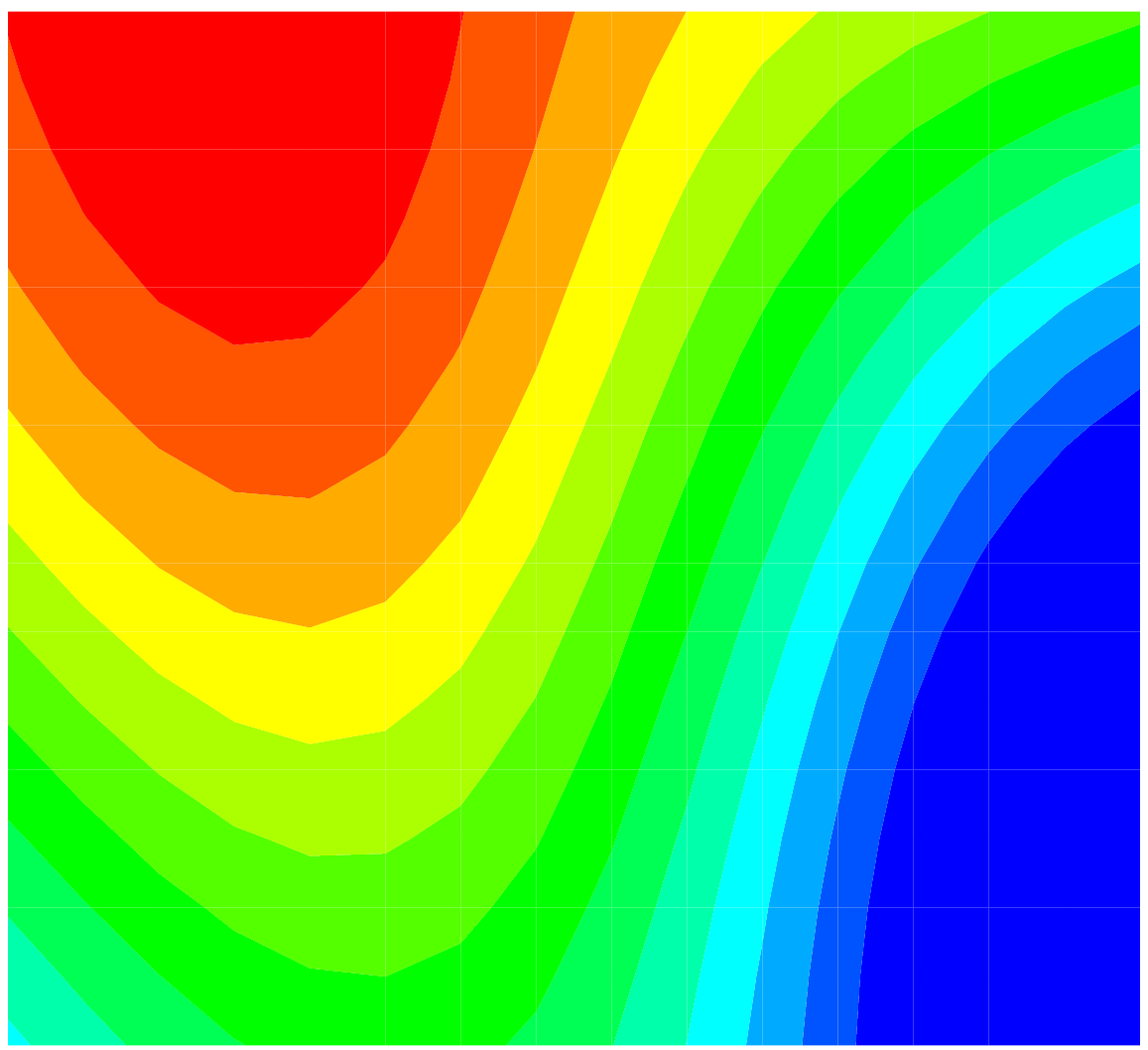}
	\includegraphics[width=1.575in]{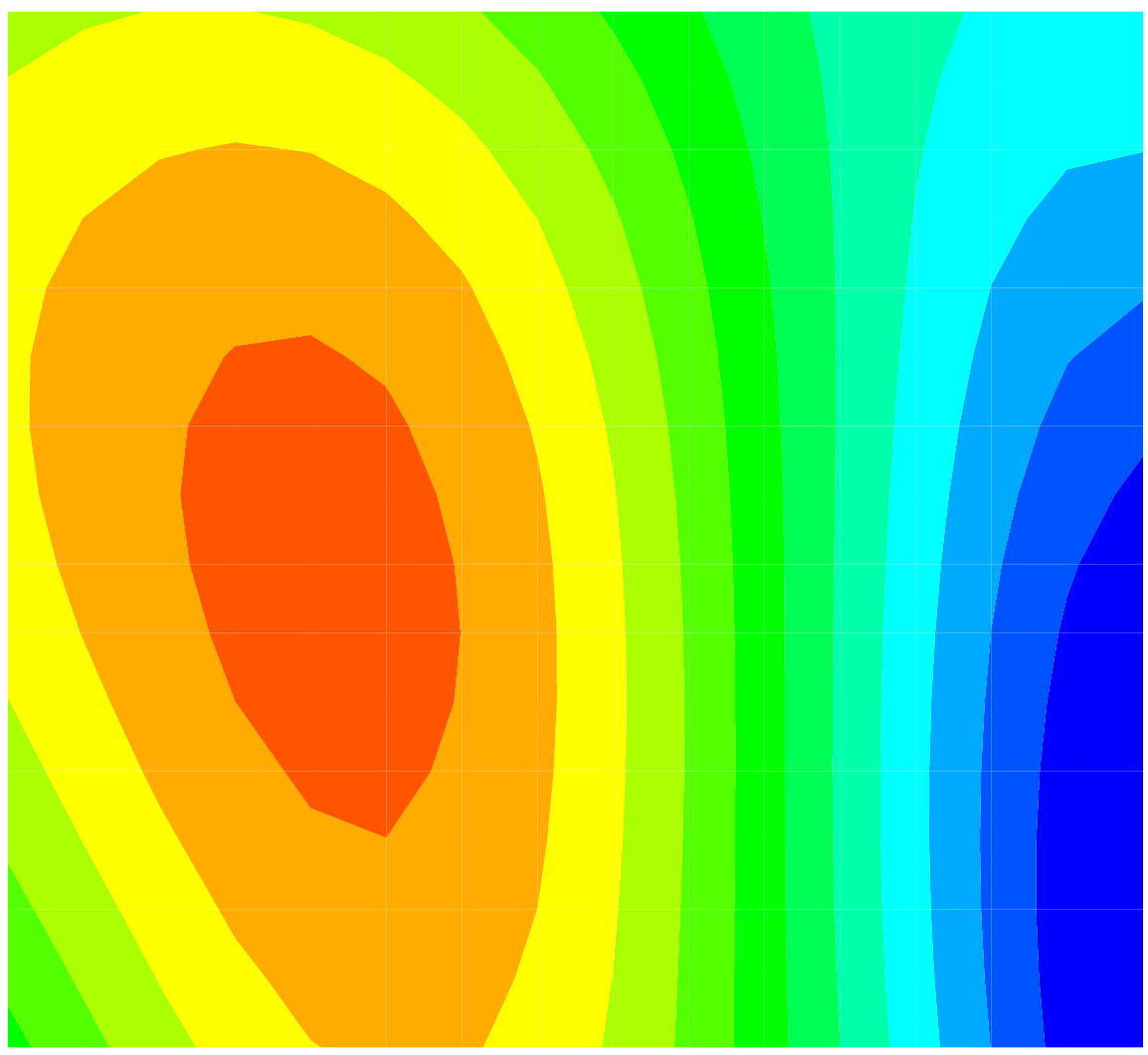}\\
	\vspace{0.3cm}
	\includegraphics[width=1.575in]{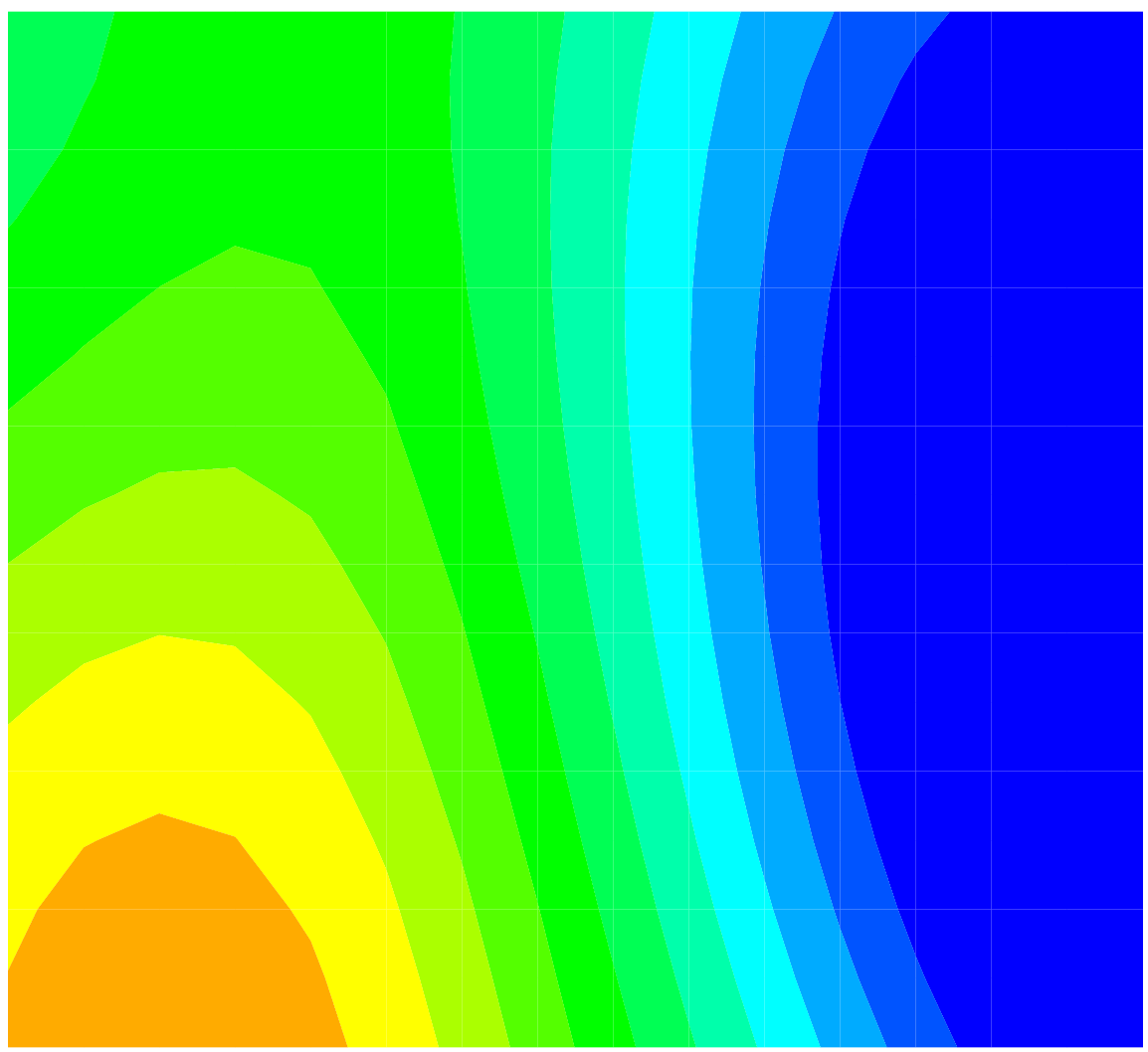}
	\includegraphics[width=1.575in]{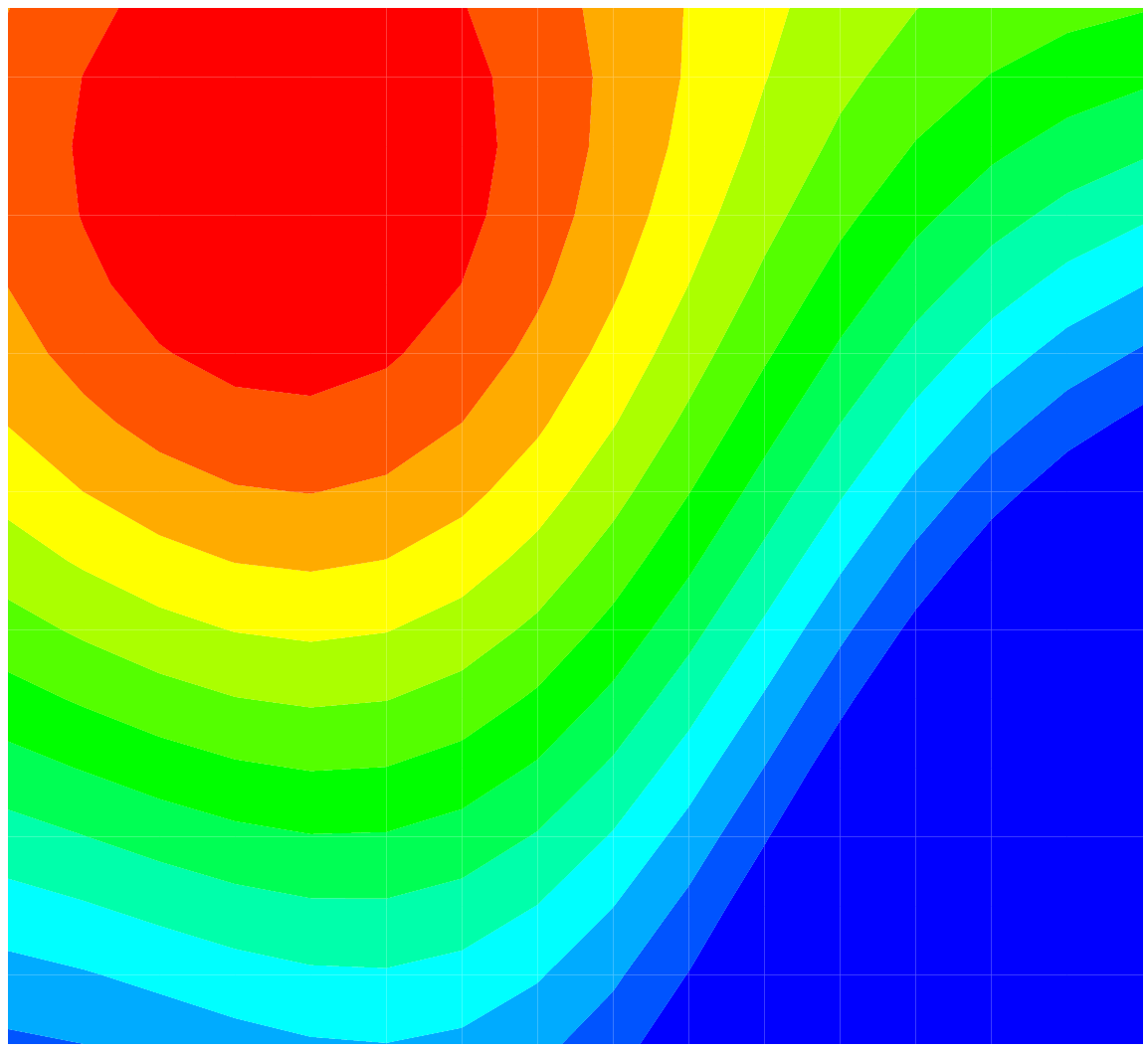}
	\includegraphics[width=1.575in]{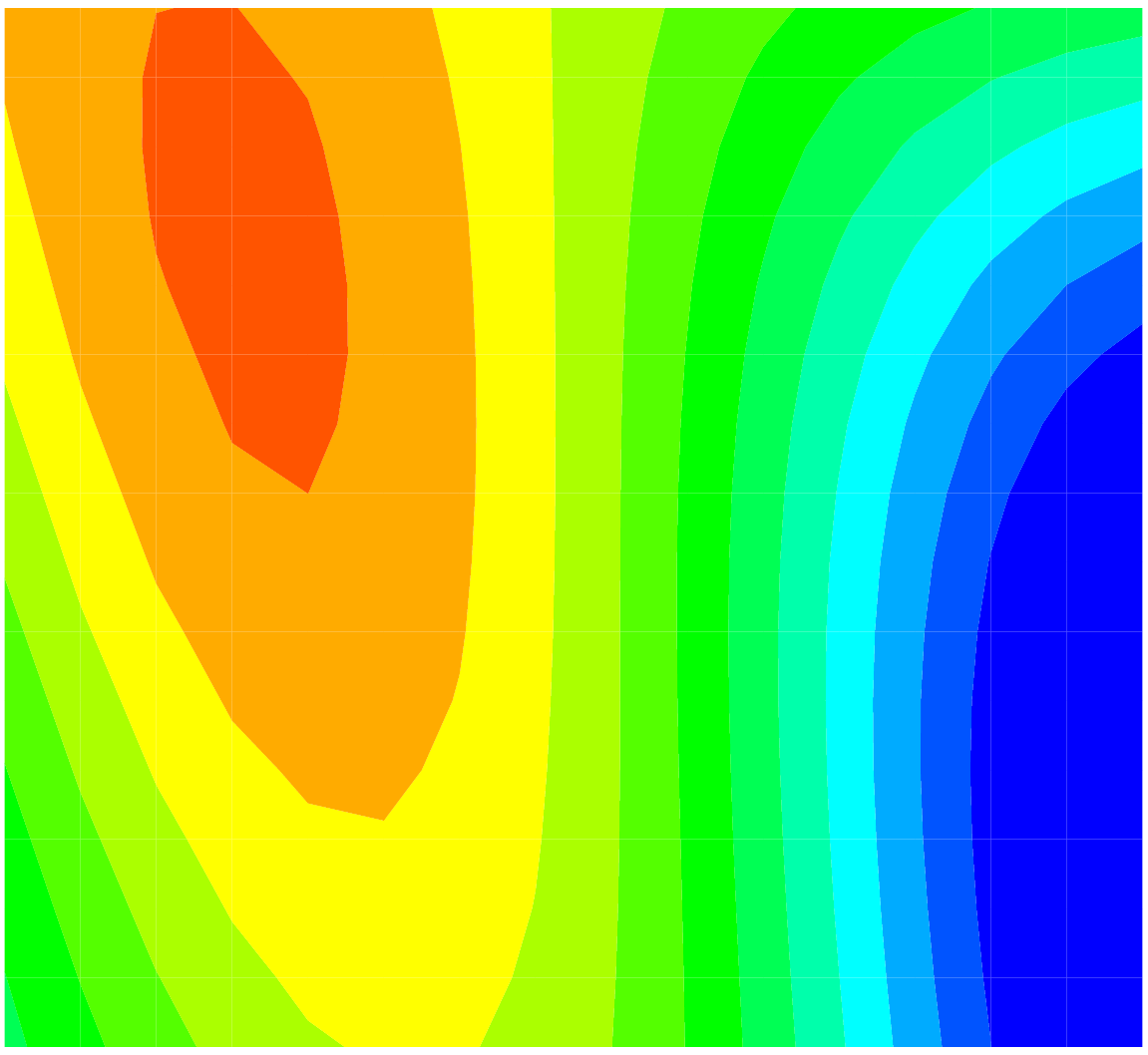}\\
	\vspace{0.3cm}
	\includegraphics[width=1.575in]{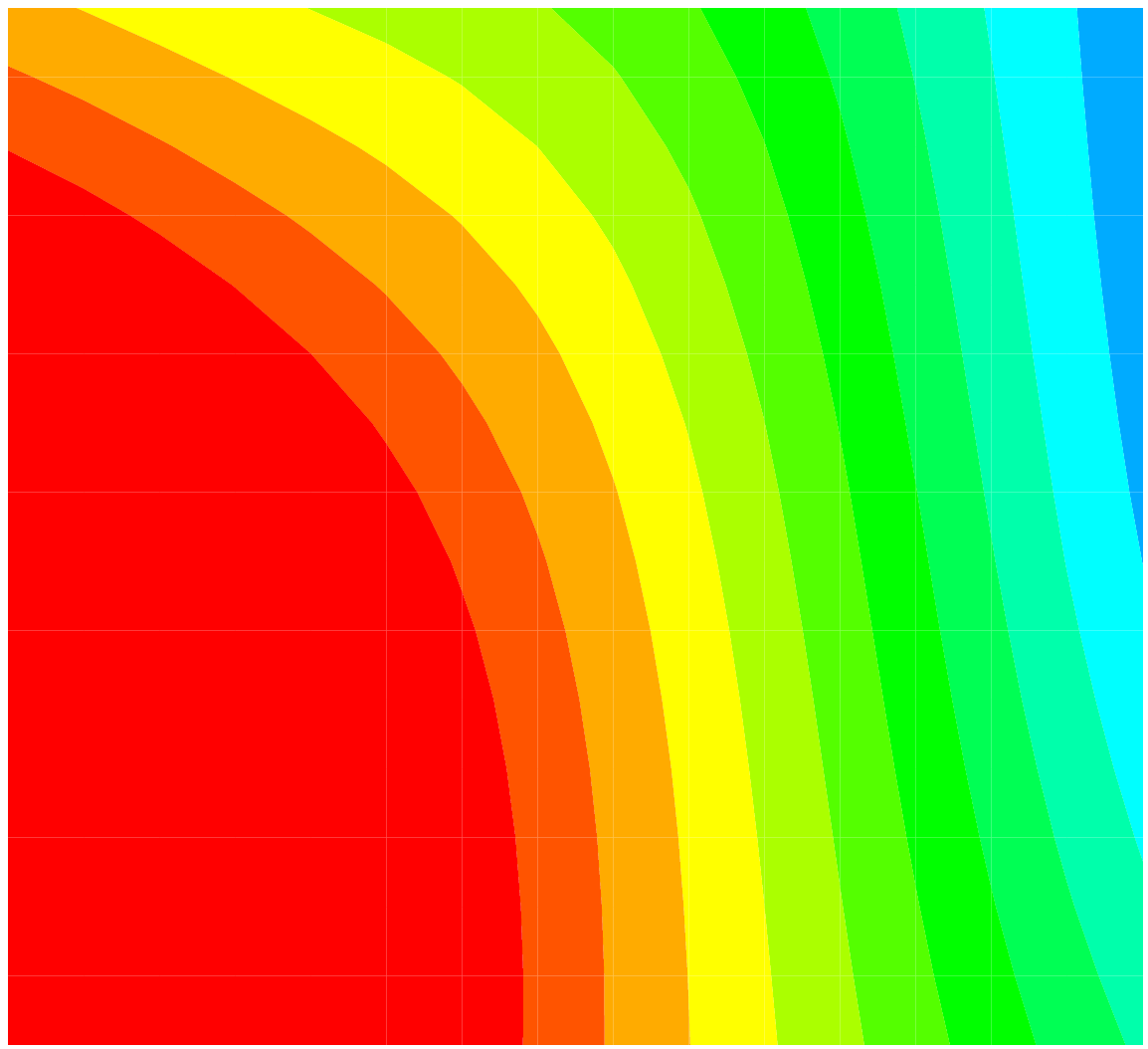}
	\includegraphics[width=1.575in]{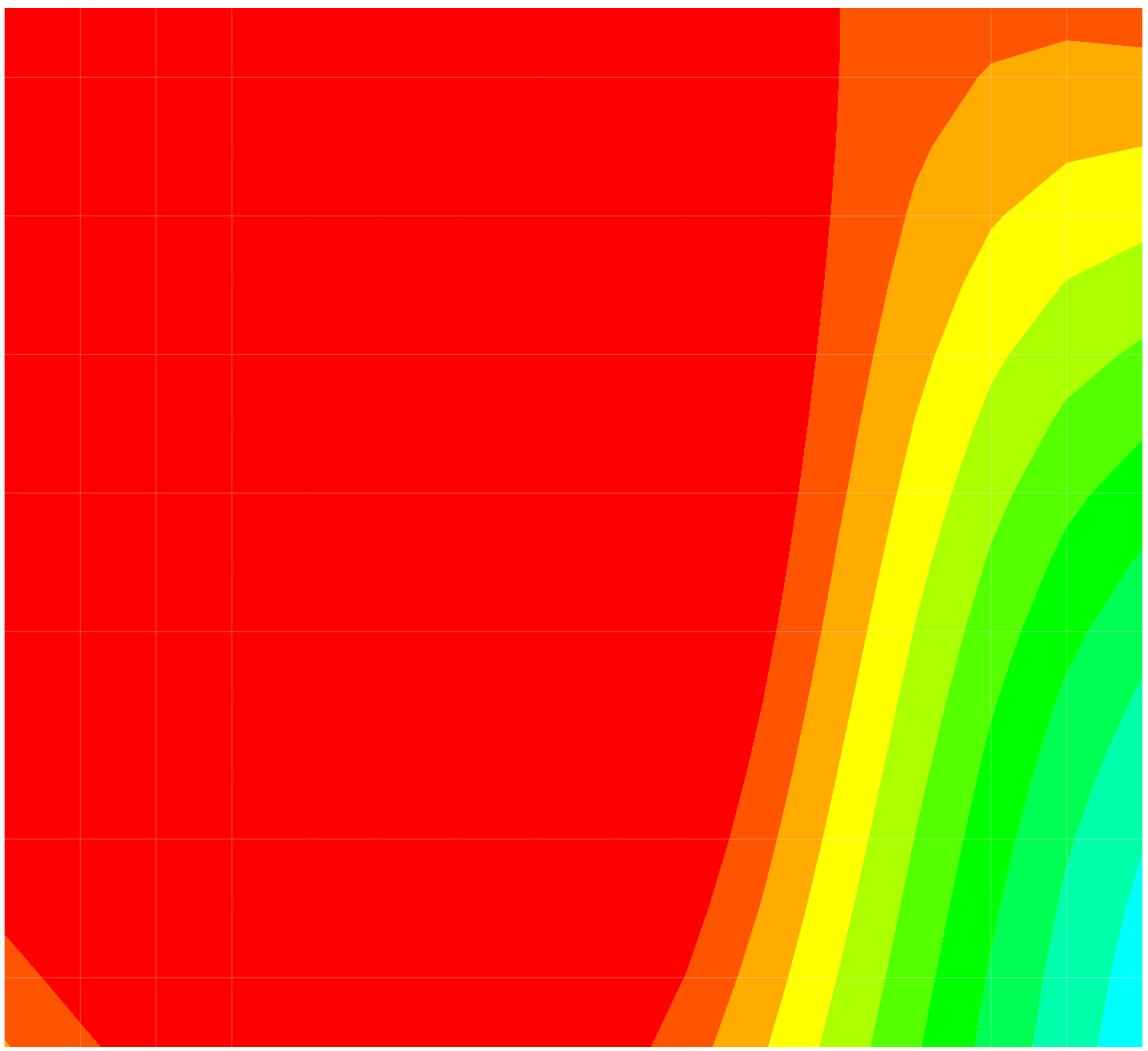}
	\includegraphics[width=1.575in]{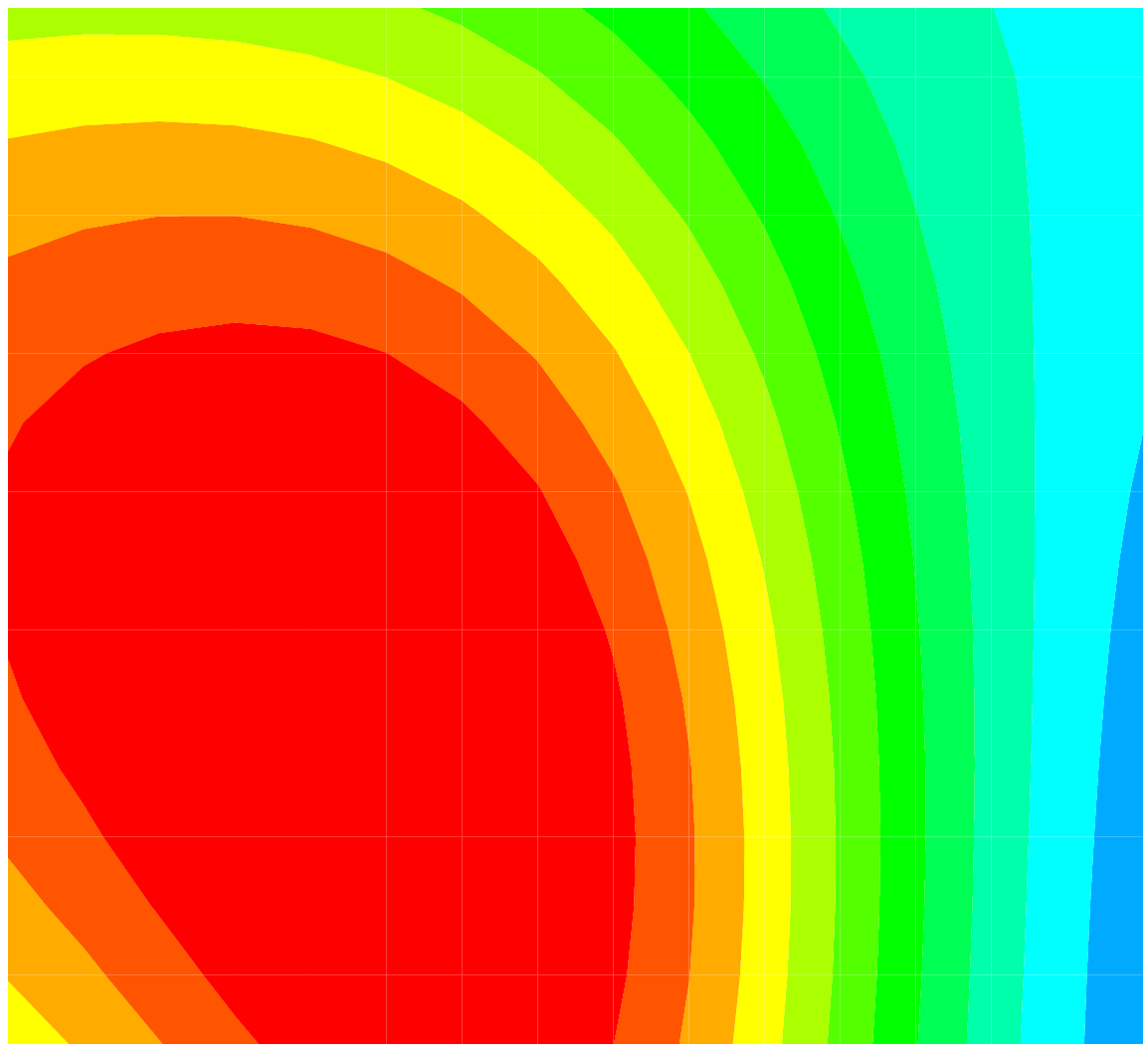}\\ 
	\vspace{0.3cm}
	\includegraphics[width=1.575in]{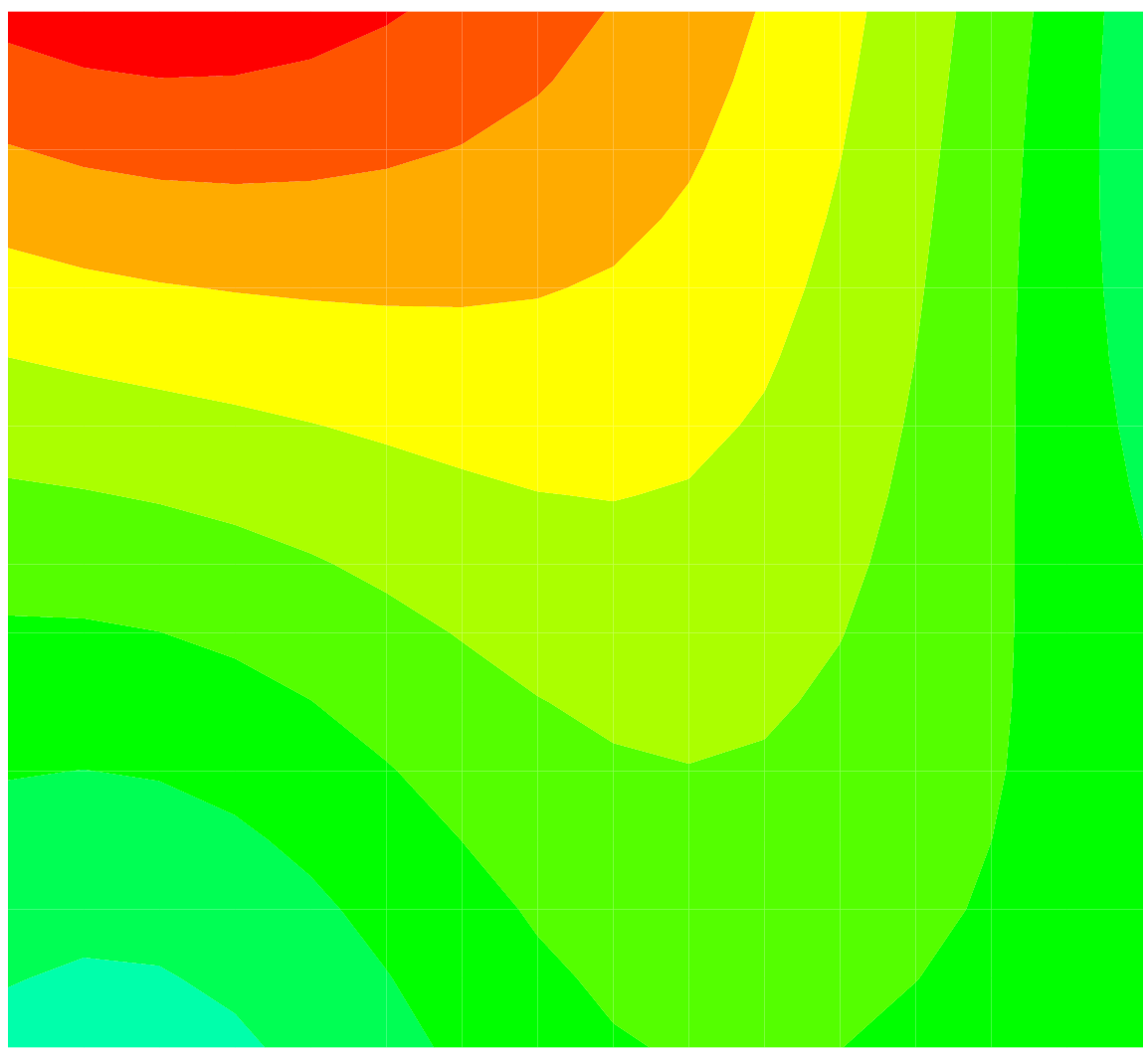}
	\includegraphics[width=1.575in]{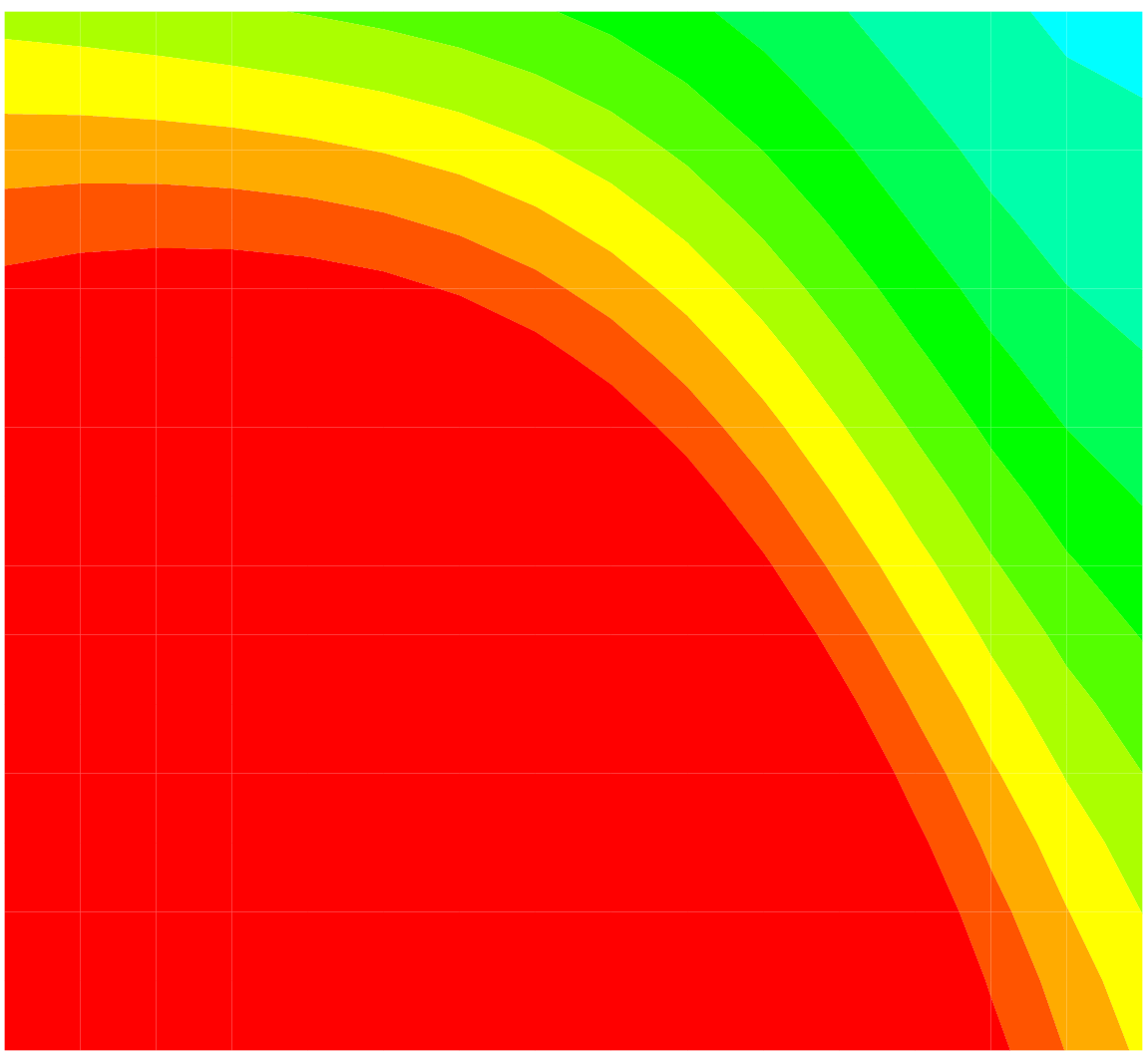}
	\includegraphics[width=1.575in] {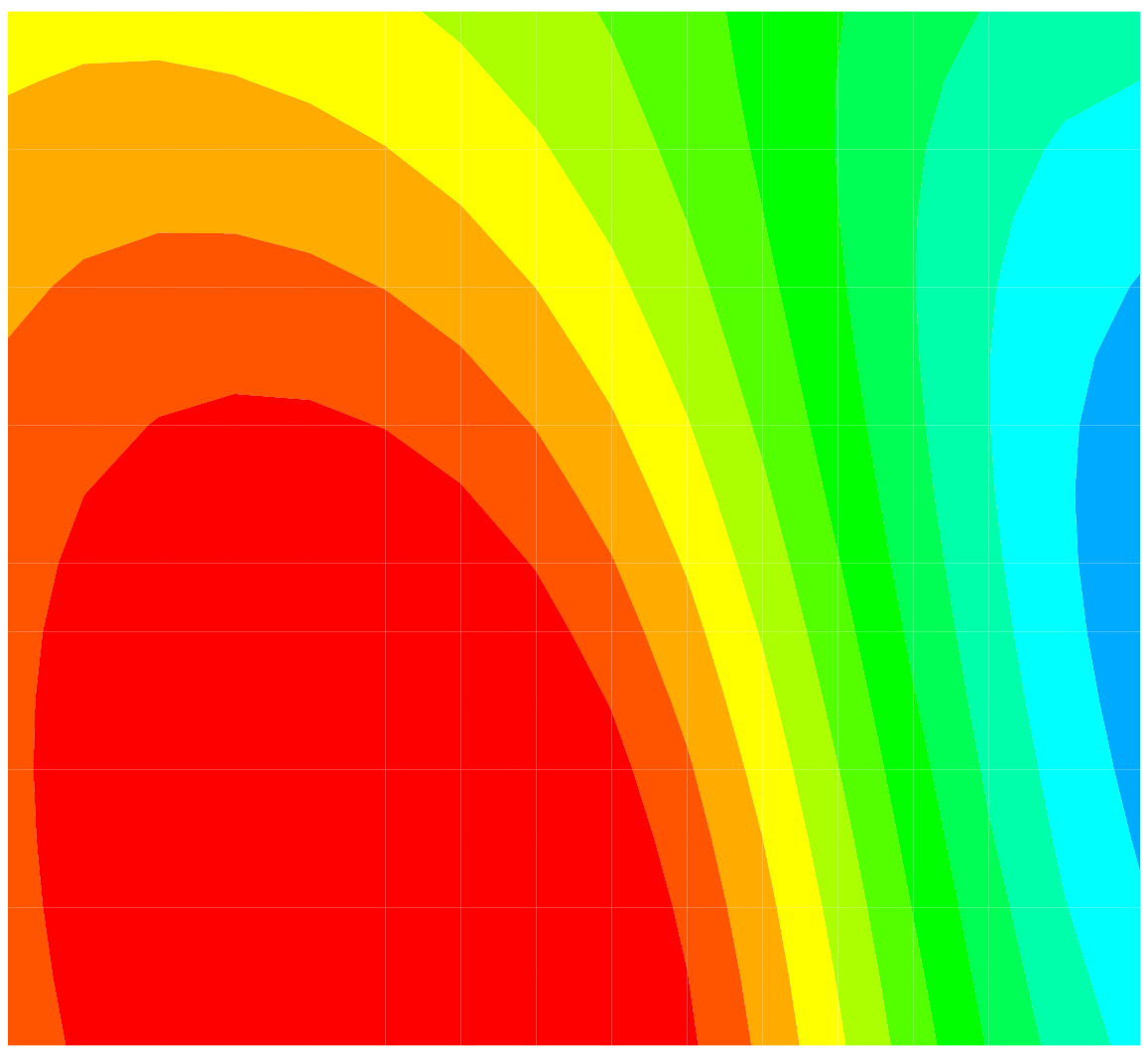}       		
	\caption{ Log permeability fields recovered in the MCMC without conditioning. Rows refer to four consecutive chains and columns refer to iterations 40, 5000 and 10000, respectively.}
	\label{perms-two-stage}
\end{figure}

\begin{figure}[H]
	\centering
	\includegraphics[width=1.575in]{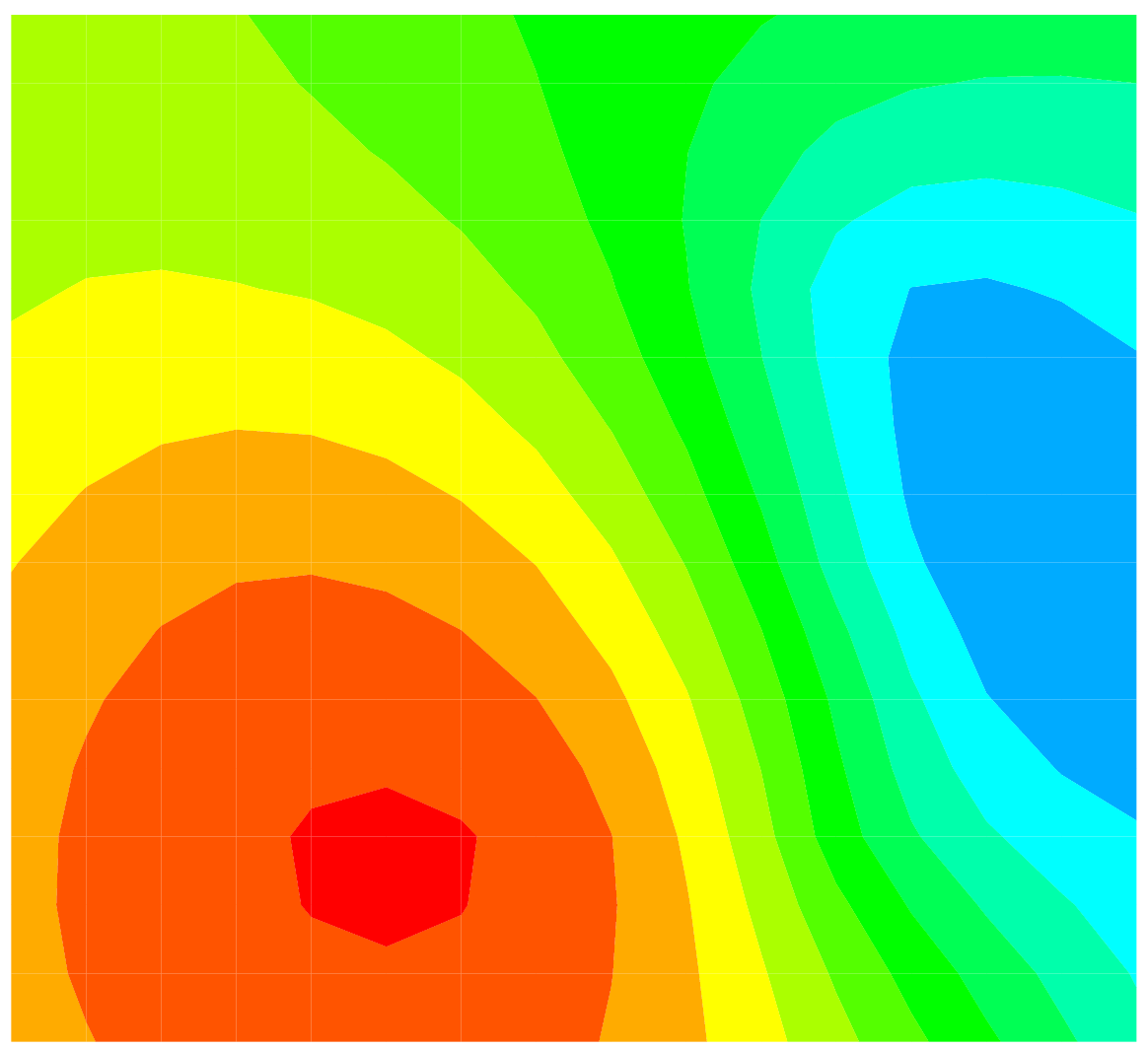}
	\includegraphics[width=1.575in]{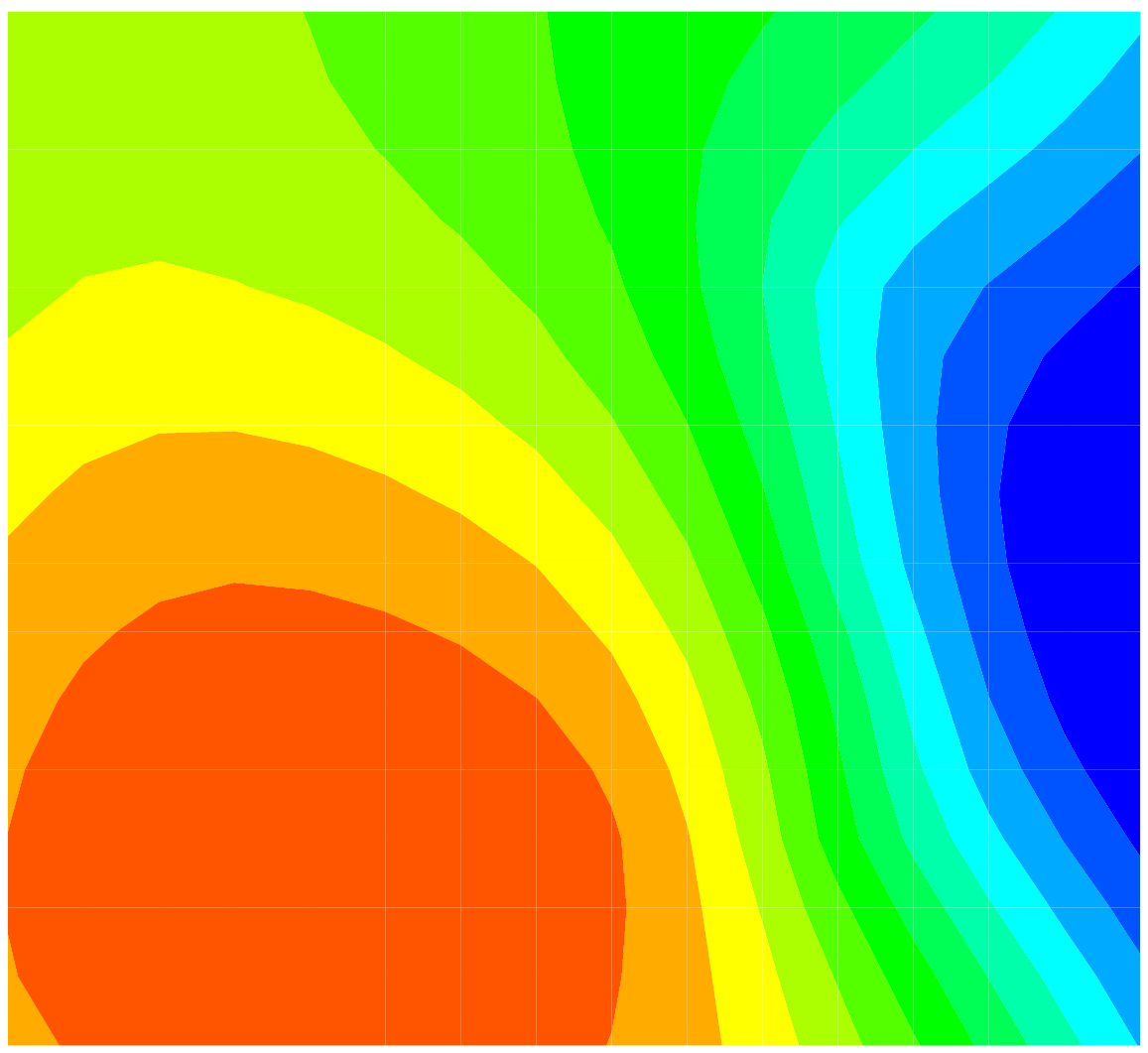}
	\includegraphics[width=1.575in]{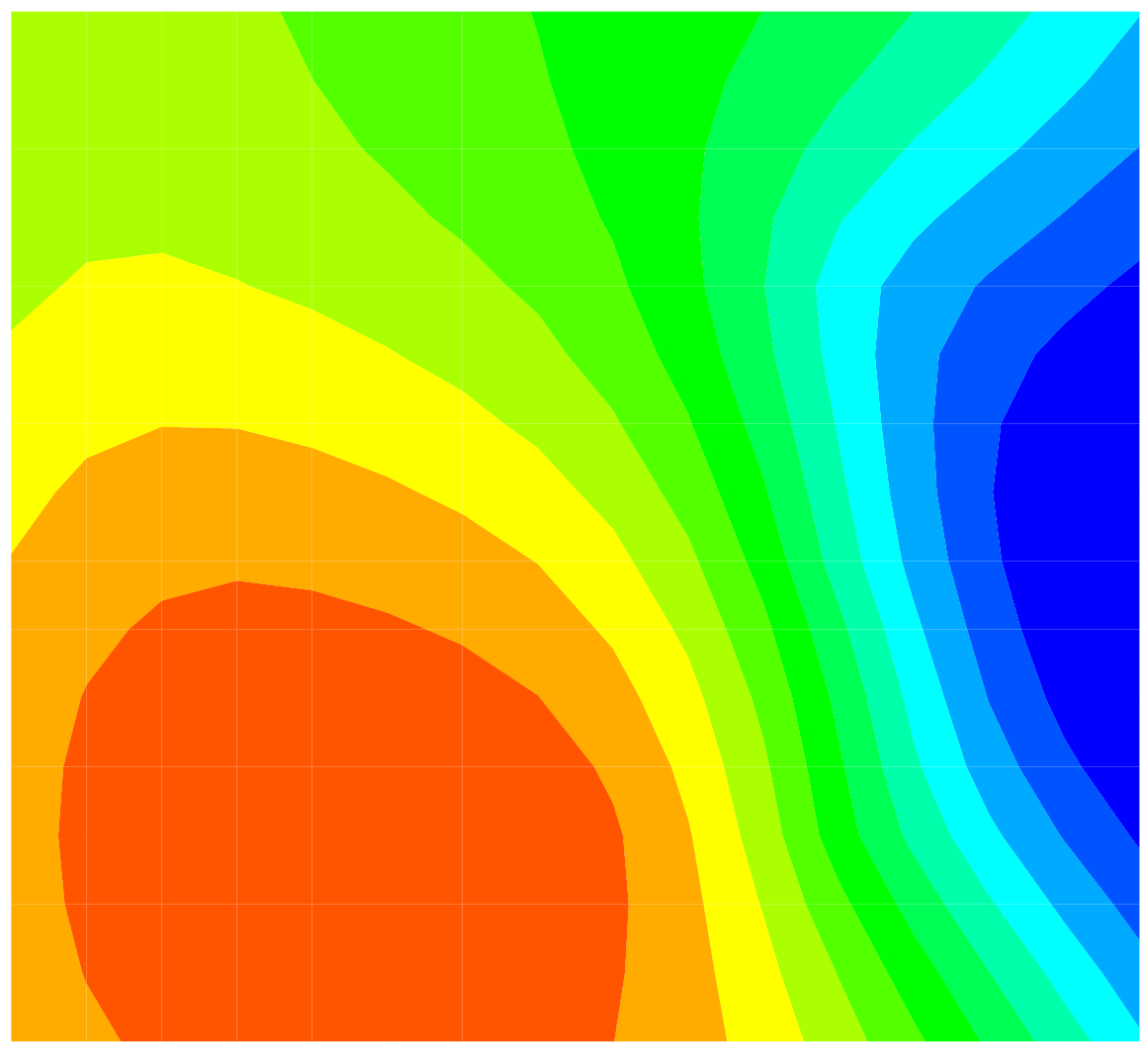}\\
	\vspace{0.3cm}
	\includegraphics[width=1.575in]{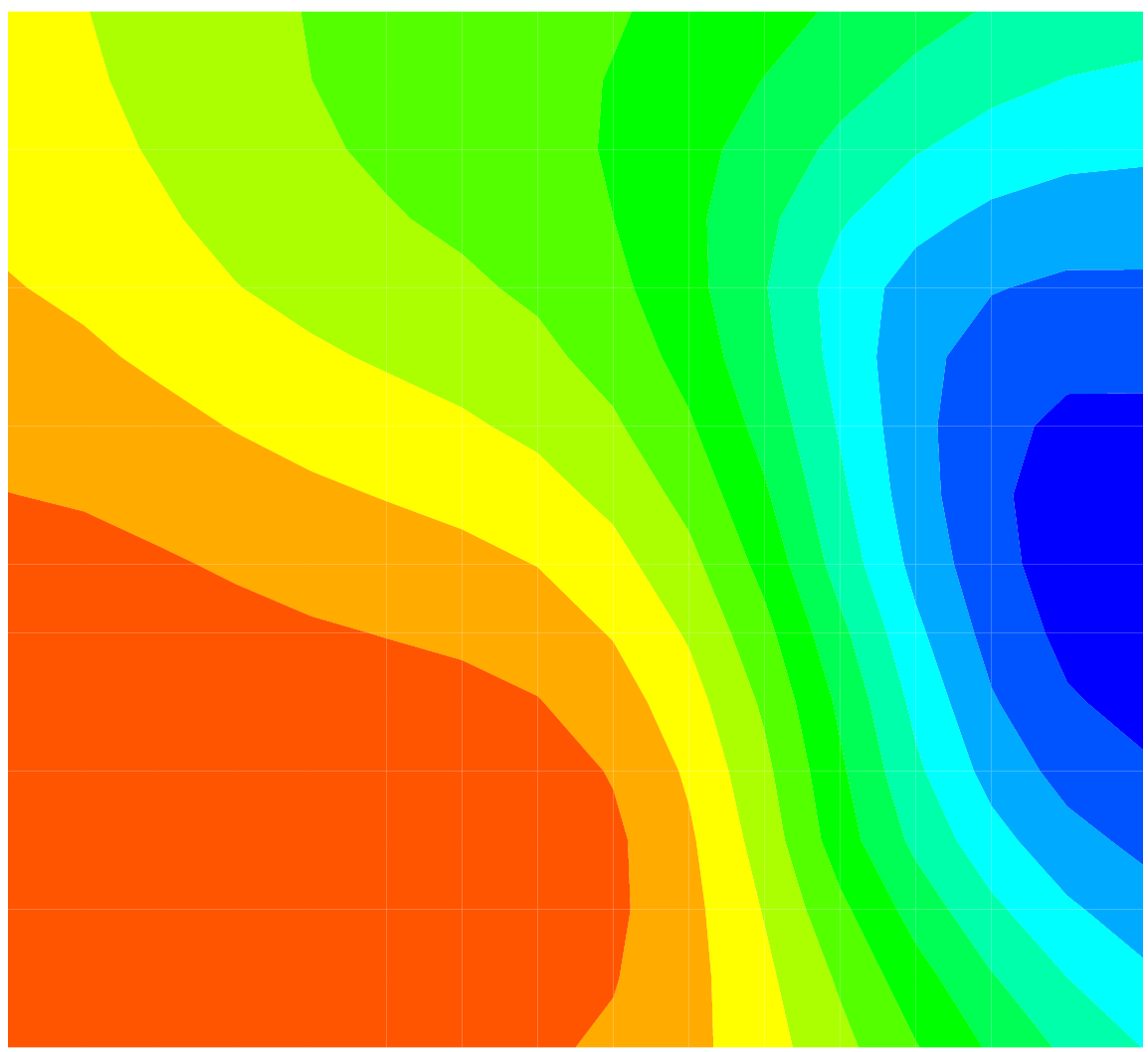}
	\includegraphics[width=1.575in]{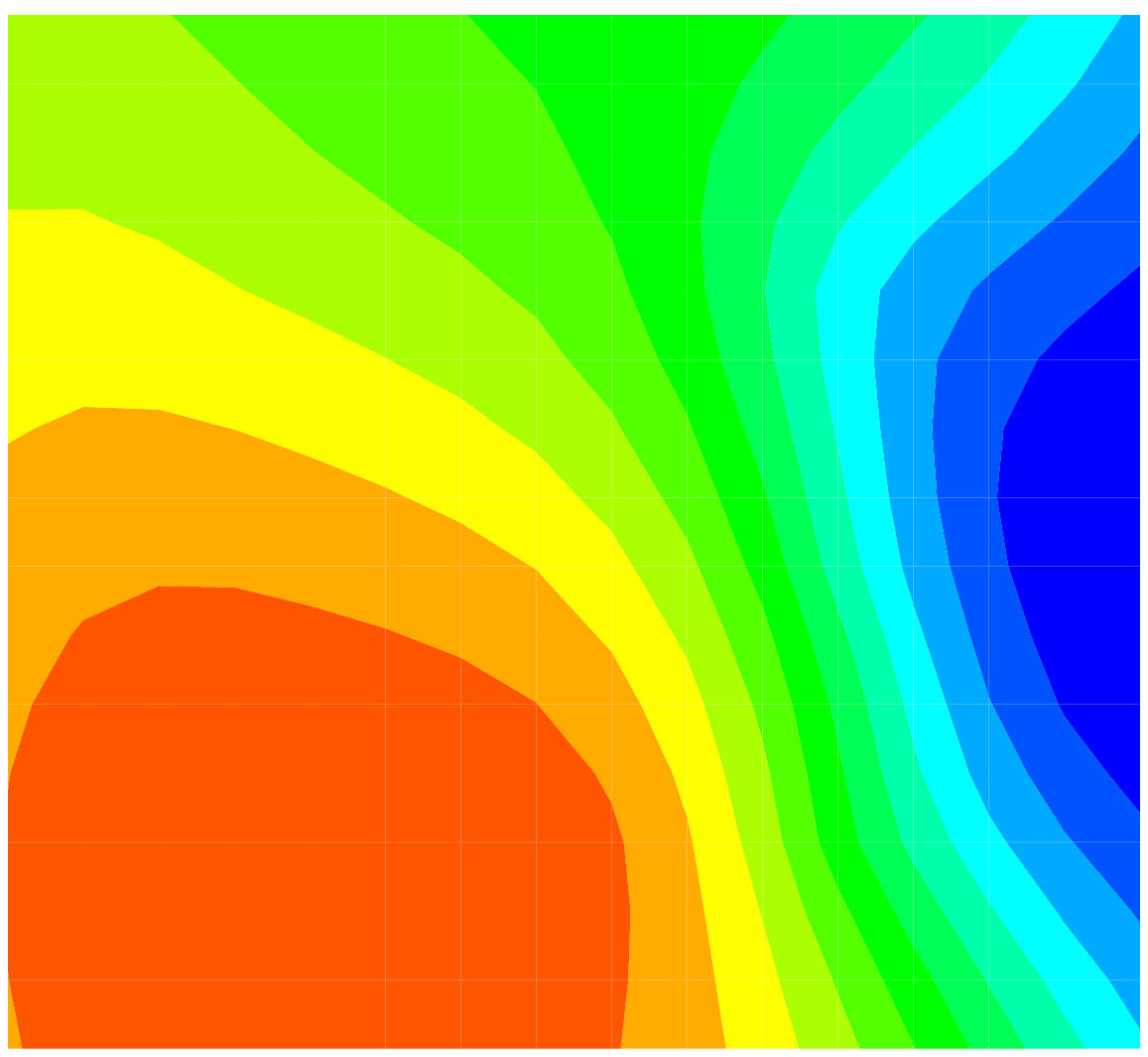}
	\includegraphics[width=1.575in]{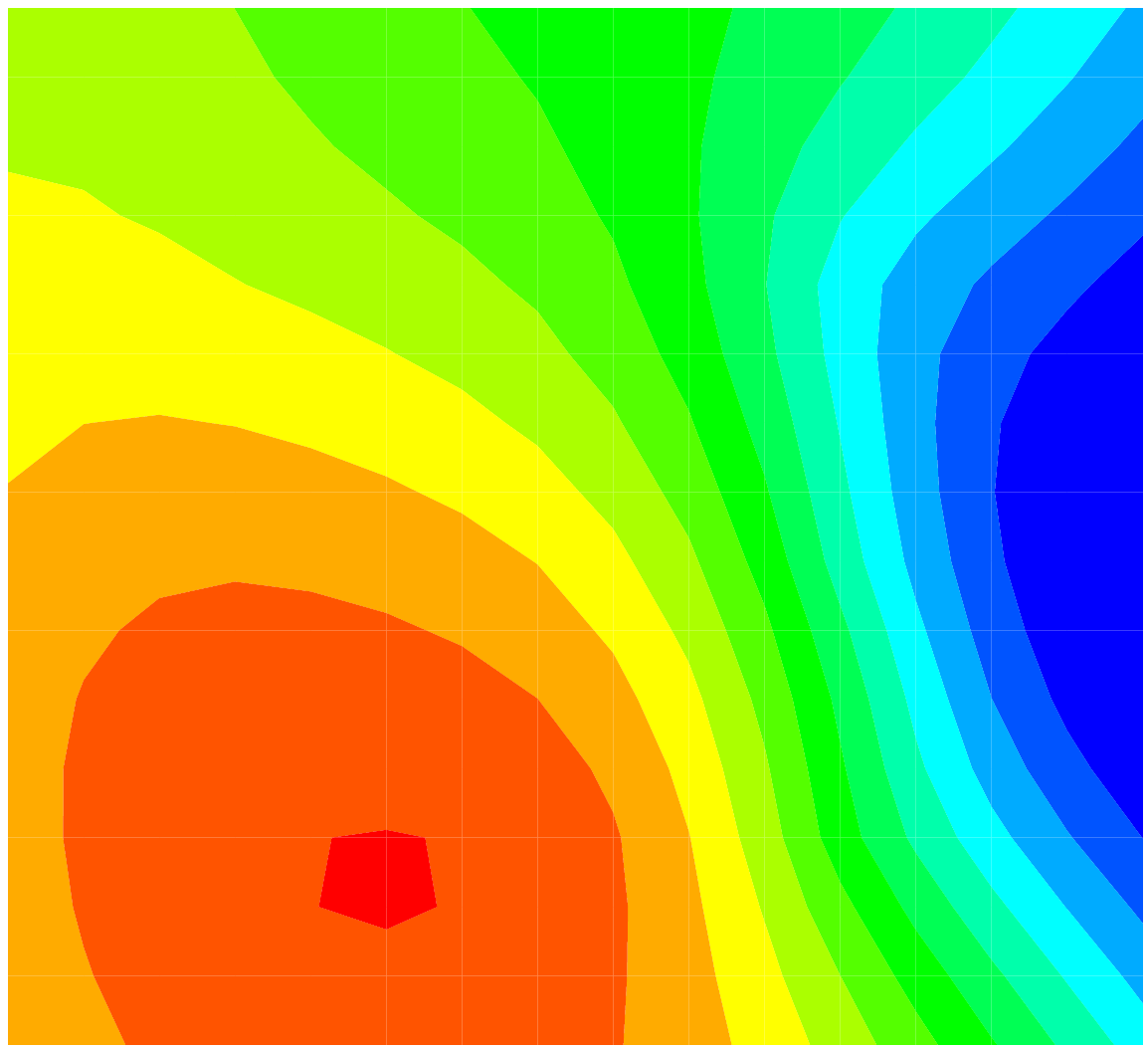}\\
	\vspace{0.3cm }
	\includegraphics[width=1.575in]{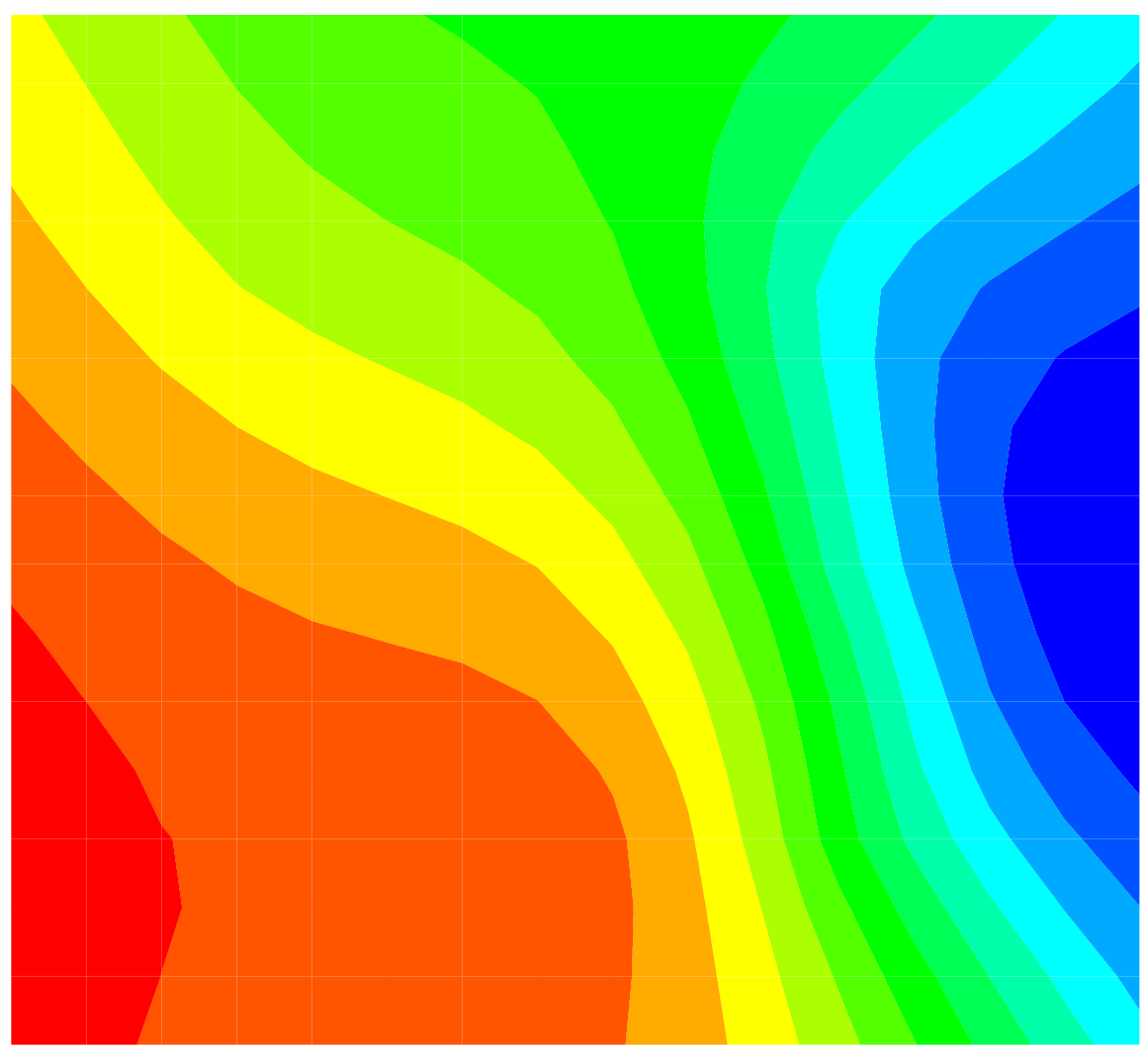}
	\includegraphics[width=1.575in]{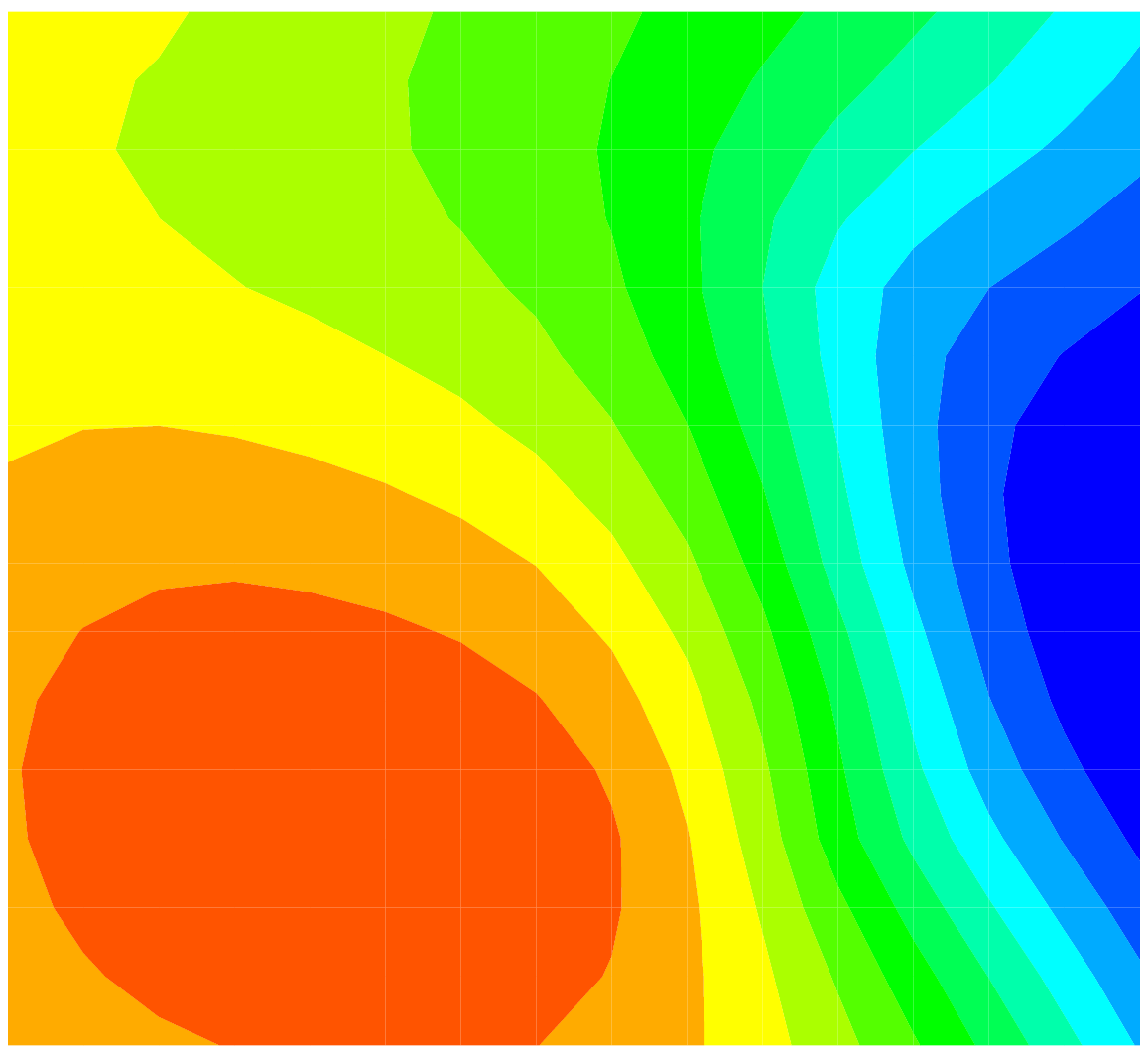}
	\includegraphics[width=1.575in]{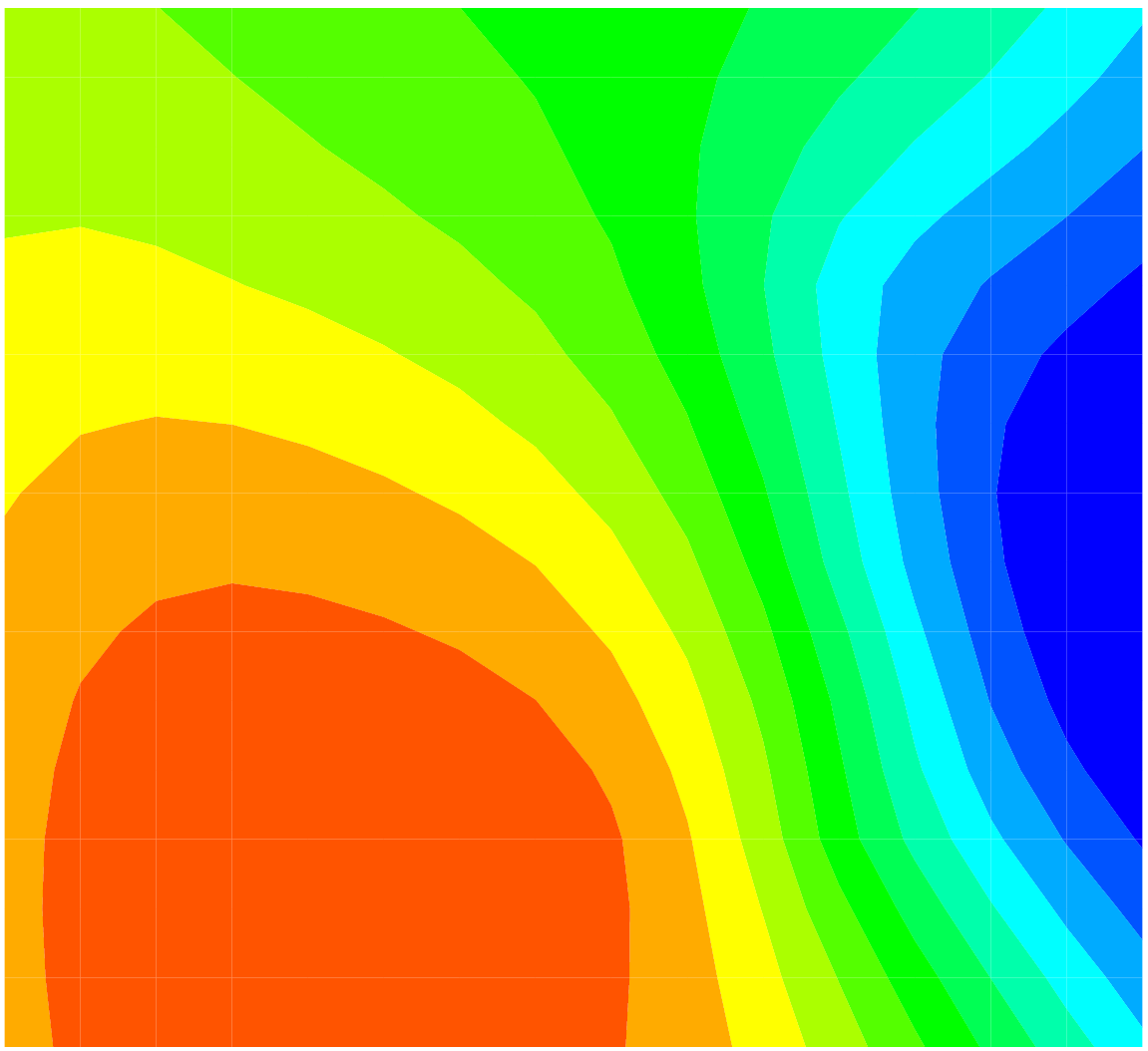}\\
	\vspace{0.3cm }     
	\includegraphics[width=1.575in]{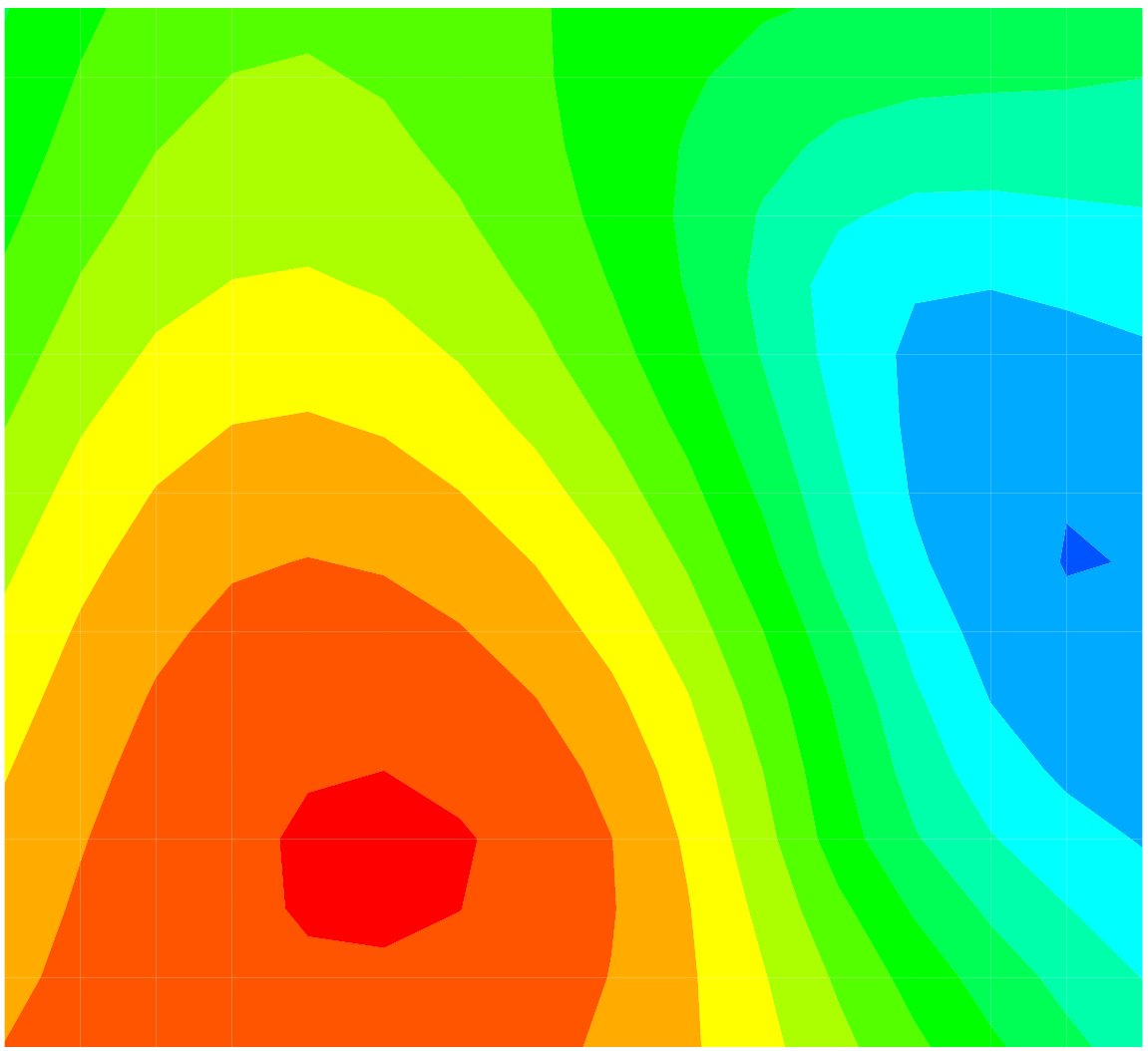}
	\includegraphics[width=1.575in]{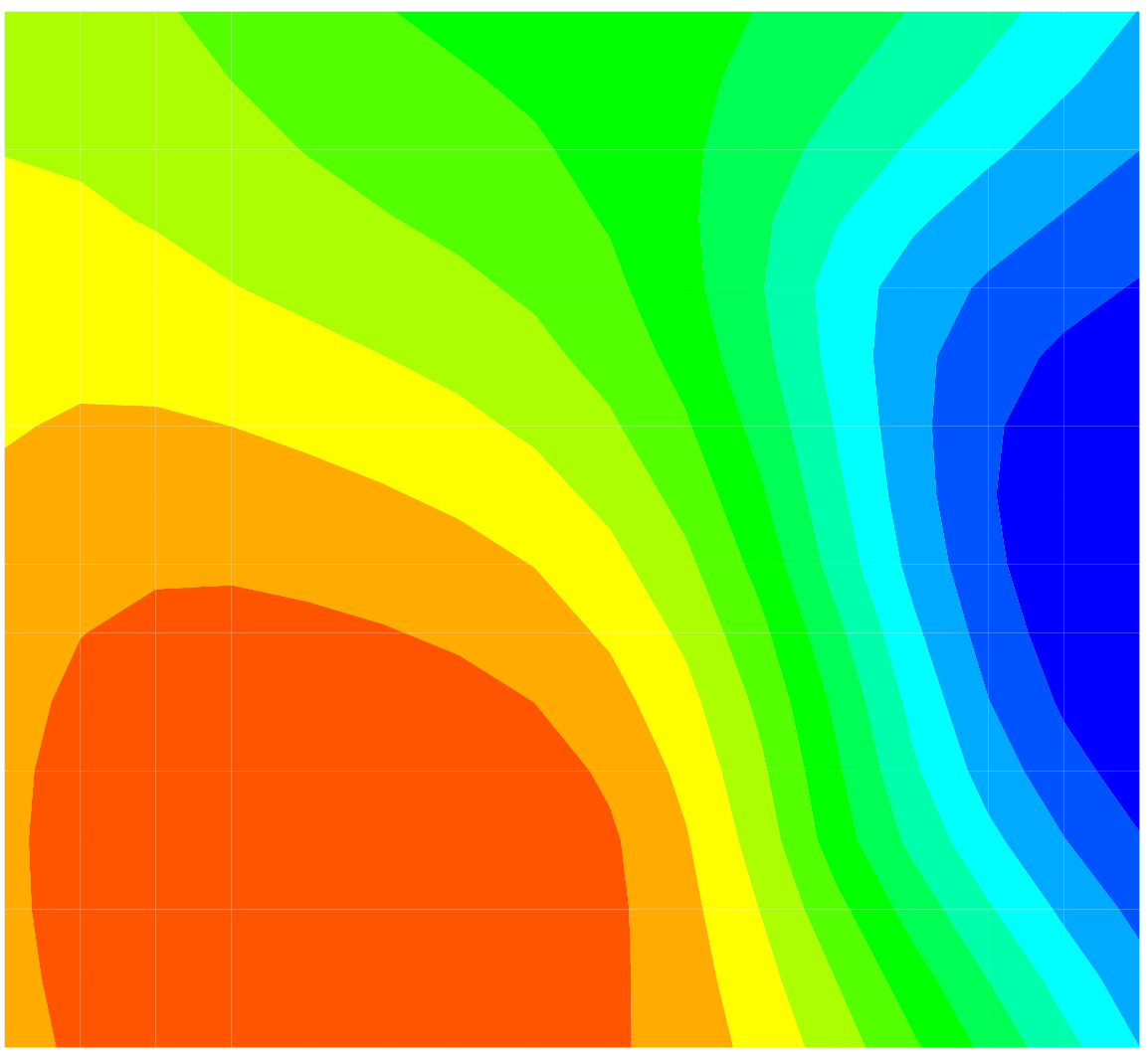}
	\includegraphics[width=1.575in] {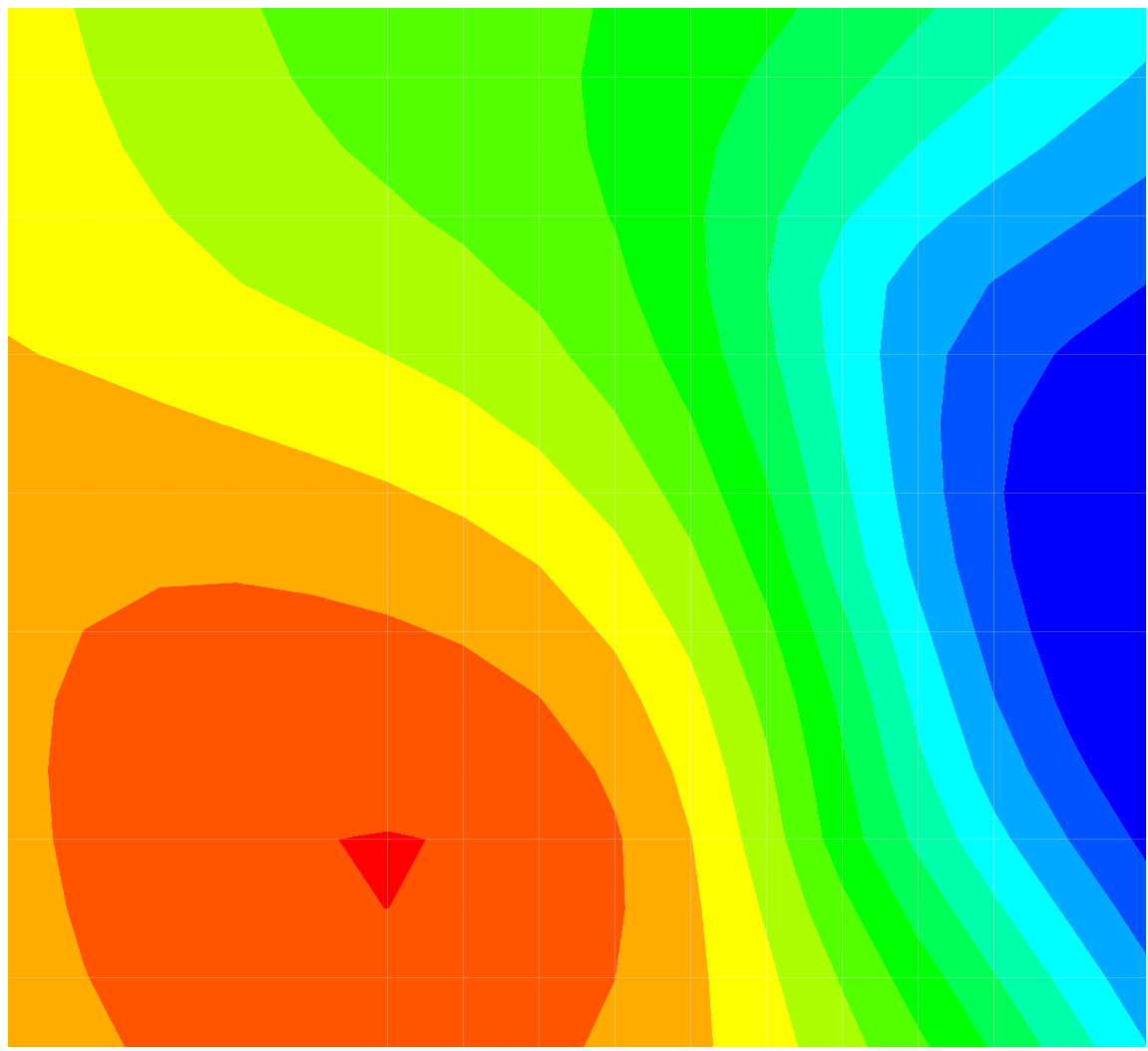}
	\caption{ Log permeability fields recovered in the MCMC with conditioning. Rows refer to four consecutive chains and columns refer to iterations 40, 5000 and 10000, respectively.}
	\label{perms_conditioning}
\end{figure}

Now, we discuss the permeability fields that are accepted in the MCMC methods. In Fig.~\ref{perms-two-stage} we show twelve permeability fields (three per chain) at 40, 5000, and 10000 iterations, respectively, in four chains for the MCMC method without conditioning. In Fig.~\ref{perms_conditioning} we present the corresponding permeability fields in four chains obtained from the MCMC method with conditioning. The process of constructing the permeability fields for this MCMC method with conditioning has two steps: In the first step, we create a permeability field using the KLE in which we provide the projected proposal as an input. Then we add the kriging data obtained from the measurements at nine locations (See section~\ref{proj}). We observe that the MCMC method with conditioning recovers the reference permeability field in Fig.~\ref{ref_perm} much better than the MCMC method without conditioning. This result is consistent with the MCMC convergence discussed in the above paragraph. Thus, we conclude that the MCMC method with conditioning shows a better performance for the characterization of the permeability field in our problem.

\section{Conclusions}
\label{conclusions}
In this work, we considered the use of MCMC methods for the characterizations of subsurface formations. Usually, sparse measurements of quantities of interest (such as the permeability field) are available in several field locations and they have to be incorporated in the characterization procedure to reduce uncertainty as much as possible.

We introduced a computationally efficient method for conditioning (the log of) permeability fields. The method is based on a projection onto the nullspace of a data matrix defined in terms of a KL expansion. We also presented multi-chain studies that illustrate the importance of conditioning in accelerating the MCMC convergence. 

The authors and their collaborators are currently applying the Multiscale Perturbation Method ~\cite{mpm_2020,mpm_2021} to further speed up MCMC studies.

\section*{Acknowledgment}
A. Rahunanthan was supported by grants from the National Science Foundation (Grant No. HRD--1600818), and NIFA/USDA through Central State University's Evans-Allen Research Program (Fund No. NI201445XXXXG018-0001).

All the numerical simulations presented in this paper were performed on the GPU Computing cluster housed in the Department of Mathematics and Computer Science at Central State University.
%
%
%

\bibliographystyle{splncs04}
\bibliography{porous_bib,REFERENCIAS}

\end{document}